\documentclass[final,3p,times,sort&compress,doublespacing,fleqn]{elsarticle}

\makeatletter
\def\ps@pprintTitle{%
  \let\@oddhead\@empty
  \let\@evenhead\@empty
  \def\@oddfoot{\reset@font\hfil\thepage\hfil}
  \let\@evenfoot\@oddfoot
}
\makeatother

\usepackage{multirow,multicol,setspace}

\usepackage{amsthm,amssymb,array,float,latexsym,amsmath,amssymb,amsbsy,eurosym,theorem,amsfonts,mathrsfs,bbm,verbatim,dsfont}
\usepackage{booktabs}
\usepackage{caption}
\usepackage{subcaption}

\usepackage{rotating}
\usepackage{tabularx}
\usepackage[table]{xcolor}
\usepackage[]{graphicx}
\usepackage{color}

\usepackage[most]{tcolorbox}
\usepackage{titlesec}

\usepackage{scalerel,stackengine}

\titleformat{\paragraph}
{\normalfont\normalsize}{\theparagraph}{1em}{}
\titlespacing*{\paragraph}
{0pt}{3.25ex plus 1ex minus .2ex}{1.5ex plus .2ex}

\renewcommand{\b}[1]{\boldsymbol{#1}} 
\renewcommand{\c}[1]{\mathcal{#1}}

\renewcommand{\v}[1]{\mathbbm{#1}} 
\renewcommand{\v}[1]{\text{\usefont{U}{bbm}{m}{n}#1}} 

\renewcommand{\o}[1]{\overline{#1}}

\newcommand{\lp}{\left(}
\newcommand{\rp}{\right)}

\newcommand{\lb}{\left[}
\newcommand{\rb}{\right]}

\newcommand{\p}{\partial}

\newcommand{\f}{\displaystyle\frac}

\newcommand{\bnull}{\b0}

\renewcommand{\d}{\mbox{d}}

\newcommand{\grad}{\ensuremath{\mbox{grad}}}

\newcommand{\dyad}{\otimes}

\newcommand{\assembly}[2]{ \mathop{\mbox{\Large {\textsf {\textbf A}}}}_{#1}^{#2} }

\newcommand{\sth}{\{\bullet\}}

\definecolor{gray}{rgb}{0.75, 0.75, 0.75}
\definecolor{yellow}{rgb}{1, 0.7, 0.2}
\definecolor{green}{rgb}{0.3, 0.9, 0.3}
\definecolor{brown}{rgb}{0.6, 0.3, 0.2}
\definecolor{magenta}{rgb}{0.9, 0.1, 0.9}
\definecolor{light}{rgb}{1, 0.7, 0.7}

\graphicspath{{figs/}}

\begin{document}

\begin{frontmatter}

\title{On continuum modeling of cell aggregation phenomena}

\author[erlangen]{Soheil~Firooz\corref{cor}}
\ead{soheil.firooz@fau.de}
\author[siemens]{Stefan~Kaessmair}
\author[biology,maxplankphysic]{Vasily~Zaburdaev}
\author[bilkent]{Ali~Javili}
\author[erlangen,glasgow]{Paul~Steinmann}
\address[erlangen]{Institute of Applied Mechanics, University of Erlangen-Nuremberg, Egerland Str. 5, 91058 Erlangen, Germany}
\address[siemens]{Siemens Industry Software GmbH, Nordostpark 3, 90411 Nuremberg, Germany}
\address[biology]{Department of Biology, University of Erlangen-Nuremberg, 91058 Erlangen, Germany}
\address[maxplankphysic]{Max Planck Zentrum f{\"u}r Physik und Medizin, 91058 Erlangen, Germany}
\address[bilkent]{Department of Mechanical Engineering, Bilkent University, 06800 Ankara, Turkey}
\address[glasgow]{Glasgow Computational Engineering Center, James Watt School of Engineering, University of Glasgow, Glasgow G12 8QQ, United Kingdom}
\cortext[cor]{Corresponding author.}

\begin{abstract}
Cellular aggregates play a significant role in the evolution of biological systems such as tumor growth, tissue spreading, wound healing, and biofilm formation.
Analysis of such biological systems, in principle, includes examining the interplay of cell-cell interactions together with the cell-matrix interaction.
These two interaction types mainly drive the dynamics of cellular aggregates which is intrinsically out of equilibrium.
Here we propose a non-linear continuum mechanics formulation and the corresponding finite element simulation framework to model the physics of cellular aggregate formation.
As an example, we focus in particular on the process of bacterial colony formation as recently studied by Kuan et al.~\cite{Kuan2021}.
Thereby we describe the aggregation process as an active phase separation phenomenon.
We develop a Lagrangian continuum description of the problem which yields a substantial simplification to the formulations of the governing equations.
Due to the presence of spatial Hessian and Laplacian operators, a gradient-enhanced approach is required to incorporate $\c{C}^{1}$ continuity.
In addition, a robust and efficient finite element formulation of the problem is provided.
Taylor--Hood finite elements are utilized for the implementation to avoid instabilities related to the LBB condition.
Finally, through a set of numerical examples, the influence of various parameters on the dynamics of the cellular aggregate formation is investigated.
Our proposed methodology furnishes a general framework for the investigation of the rheology and non-equilibrium dynamics of cellular aggregates.
\end{abstract}


\end{frontmatter}

\section{Introduction}\label{sec:intro}
\noindent
Within the human body, most cells interact with their neighboring cells and with their extracellular matrix to establish a unique organization.
These cell-cell and cell-matrix interactions form a complex network of mechanical, biological and chemical signals which play a significant role in cell physiology~\cite{Futaki2021,Griffith2006,Kim2020,Pampaloni2007}.
The study of cellular interactions provides a significant insight towards a better understanding of many biological processes such as tumor growth~\cite{Kleiser1987,Eguchi2018}, tissue spreading~\cite{Douezan2012,Beaune2014,Beaune2017,Douezan2011}, biofilm formation~\cite{Clevers2016,Johnson2008,Ben-Jacob1998} and wound healing~\cite{Armstrong1992,Bi2016}.
Although there exist numerous contributions on experimental methods to examine cellular interactions and their final outcomes~\cite{Freitag1995,Brossay1994,Todd1984,Klausen2003,Imhaus2014}, there are still certain details and information that remain hard to assess.
Examples of such details include the exact relation between the interactions and individual cell's adhesion properties, or the speed at which sorting between the cells occurs.
Clearly, determining the relationships between all the problem variables via conducting experiments is an arduous and time consuming task.
Theoretical approaches, however, provide a valuable alternative as they do not suffer from such limitations and facilitate parametric studies.
Mathematical modeling allows us to quantify parameters such as local speeds, fluxes and their rates of change in a natural way.
They also enable us to investigate the effects of interactions between cells of different types or to allow the cell properties to vary independently of one another.
Given these benefits, mathematical modeling proves to be an inclusive and efficient alternative for analyses of cellular interactions~\cite{Armstrong2006}.

There exist two major approaches for mathematical description of cellular interactions: \emph{agent-based} and \emph{continuum models}.
In agent-based approaches, the cells are modeled and treated individually using a set of biophysical rules~\cite{Drasdo2005,Griffe2020,Anderson2000,Block2007}.
Each cell is approximated, for example, as a homogeneous, isotropic, elastic, spherical object parameterized by measurable biophysical and biological quantities.
The agent-based approach is particularly useful when one wants to study the interaction of individual cells with each other and with their environment~\cite{MOREIRA2002,Alber2002,Alarcon2004}.
Since this method is based on a series of rules for each cell, translating biological processes into a model is straightforward~\cite{Geris2011}.
For small-scale studies or cases in which the properties of the cells vary over distances comparable to the size of a cell, a higher degree of spatial resolution is obtained via agent-based models in comparison to continuum models.
Despite all the precision that the agent-based method offers, this approach is difficult to study analytically and its computational cost greatly increases as the number of cells increases.
For instance, simulation of a tumor growth process requires systems which evolve from a single progenitor cell to $10^{6}$ cells \emph{in vitro} and $10^{11}$ cells \emph{in vivo}.
Carrying out a computational analysis on such large cell population sizes is a cumbersome task, if possible at all~\cite{Macklin2007}.
Additionally, it is often neither desirable nor necessary to track each individual cell within a very large population.
For larger-scale applications, a continuum modeling proves to be a more viable alternative.
This approach is well suited to describe large scale phenomena where the cell properties vary smoothly over a length scale of several cell diameters and therefore the cell properties can be approximated by a local average.
Continuum models frequently involve ordinary and partial differential equations which are usually in the reaction-diffusion form~\cite{Araujo2004,Roose2007,Horstmann2004,Anguige2010,Anguige2009}.
Many aspects of tumor and tissue growth have been studied using continuum models~\cite{Byrne1996,Cristini2003,Macklin2007}.
A continuum description of cell motility due to cell-cell and cell-matrix interaction was presented in~\cite{Armstrong2006,Gerisch2008,Ramis-Conde2008}.
They introduced a non-local interaction term to account for adhesion between the cells and between the cells and matrix.
Coarse grained continuum approaches such as hydrodynamic theories, have also provided a powerful tool to capture large scale emergent behaviors in active cellular systems~\cite{Callan-Jones2011a,Ranft2010,Julicher2007,Popovic2017,Prost2015a}.
Further studies on continuum modeling of cellular aggregates are available in~\cite{Kopf2013,Bao2018,Arciero2011}.
We refer to~\cite{Byrne2009} for a thorough comparison between agent-based and continuum modeling of cellular aggregates.


A prototypical biological example of cellular aggregation is the formation of bacterial microcolonies and biofilms.
One of the first steps in the process of bacterial colonization of biotic and abiotic surfaces is the formation of aggregates or colonies consisting of several thousands of cells.
Usually, these microcolonies later evolve into much more complex bacterial communities, known as biofilms~\cite{Johnson2008a,Otoole2000}.
Bacterial infections involving biofilms are far more resistant to anti-microbial treatments in many cases~\cite{Hoiby2010}.
Thus, investigation of the mechanism of bacterial microcolony formation is of immense significance in the fields of medicine and engineering.
A few well-known examples of bacterial microcolonies causing dangerous microbial infections are \emph{Pseudomonas aeruginosa}~\cite{Jin2011}, \emph{Neisseria meningitidis}~\cite{Charles-Orszag2016}, \emph{Vibrio cholerae}~\cite{Toma2002} and \emph{Neisseria gonorrhoeae}~\cite{Ponisch2017}.

In this manuscript, we focus on microcolonies of \emph{Neisseria gonorrhoeae (NG)}  bacteria.
These bacterial microcolonies are the infectious units which form on human epithelial tissue and cause gonorrhoeae, the second most common sexually transmitted disease~\cite{Jonsson1994}.
Multi-scale computational simulations have been conducted recently to study biophysical aspects of NG microcolonies~\cite{Bonazzi2018,Ponisch2017,Ponisch2018}. 
The NG bacteria, as well as many other bacteria species, use multiple long and thin retractable filaments, called type IV pili, in order to interact with the environment and with each other~\cite{Howie2005}.
A series of studies have been carried out to investigate the twitching motility of bacteria mediated by type IV pili~\cite{Mattick2002,Holz2010,Biais2010,Eriksson2015,Merz2000}.
Pili can extend from the cell body, attach to the substrate and retract.
Pilus retraction generates forces which are then translated into movement of cells.
The magnitude of the forces generated by the pilus retraction are in the range of $100-180\,\text{pN}$ which is considered as one of the strongest active molecular forces known in nature~\cite{Maier2003,Ponisch2017}.
Additionally, pili of one cell could also extend and attach to pili of other cells.
Retraction of the attached pili network attracts the cells towards each other and leads to formation of an aggregate.
These cycles of growth, attachment, detachment and retraction drive the cell motility on substrates and the aggregate formation process~\cite{Zaburdaev2014}.
Pili mediated cell-cell and cell-matrix interactions are crucial for the formation and maintenance of microcolonies~\cite{Weber2015,Todd1984,Imhaus2014,Taktikos2015}.

The main objective of this contribution is to formulate and simulate the process of cell aggregation phenomena within a nonlinear continuum mechanics framework.
We develop our framework in a Lagrangian setting which yields considerable simplification of the equations and enables implicit time integration which considerably increases the computational robustness.
In doing so, we take a prototypical example of NG bacteria and we describe the process of colony formation as an active phase separation phenomenon.
Our work is mainly based on the coarse grained approach previously developed by Kuan et al.~\cite{Kuan2021,Kuan2021b}.
While we focus on intercellular interactions, our aim is to provide a robust and efficient computational setting for generic cell-matrix interaction problems, and to develop its fully nonlinear finite element implementation.
Our proposed framework provides a versatile and reliable simulation technique that allows studying the processes of aggregate formation under high forces and strong phase-separated regimes nearing much closer to the physiologically relevant conditions.

\begin{table}[b!]
\caption{Summary of key definitions and notations.}
\label{tab:notation} 
\setstretch{1.2}
\begin{tabular}{|p{0.057\textwidth}  p{0.38\textwidth} p{0.06\textwidth} p{0.42\textwidth}|}
\toprule
$\sth$                        & an arbitrary quantity                           &$\dot{\sth}$                  & material time derivative of $\sth$\\[-2.5pt]
$\nabla_{\b{X}}\sth$          & material gradient of $\sth$                     &$\nabla_{\b{x}}\sth$          & spatial gradient of $\sth$ \\[-2.5pt]
$\nabla_{\b{X}}\cdot\sth$     & material divergence of $\sth$                   &$\nabla_{\b{x}}\cdot\sth$     & spatial divergence of $\sth$ \\[-2.5pt]
$\Delta_{\b{x}}\sth$          & spatial Laplacian of $\sth$                     &$\nabla^{2}_{\b{x}}\sth$      & spatial Hessian of $\sth$ \\[-2.5pt]
$\c{L}_{t} \sth$              & Lie time derivative of $\sth$                   &$\delta\sth$                  & variation of $\sth$      \\[-2.5pt]
$R$                           & cell radius                                     &$E$                           & cell elastic modulus  \\[-2.5pt]
$f^{\text{p}}$                & pili-pili attractive force                      &$f^{\text{s}}$                & steric repulsive force \\[-2.5pt] 
$n^{\text{p}}$                & number of bound pili pairs                      &$\xi$                         & cell-substrate friction coefficient \\[-2.5pt]
$k_{\text{on}}$               & pili binding rate                               &$k_{\text{off}}$              & pili unbinding rate \\[-2.5pt]
$l$                           & pili length                                     &$\ell_0$                      & pili average length  \\[-2.5pt] 
$\b{l}_{ij}$                  & distance vector between cells $i$ and $j$       &$\b{l}$                       & spatial velocity gradient\\[-2.5pt] 
$\b{X}$                       & material position vector                        &$\b{x}$                       & spatial position vector    \\[-2.5pt]
$c_0$                         & material cell number density                    &$c_t$                         & spatial cell number density\\[-2.5pt] 
$p_0$                         & material bound pili number density              &$p_t$                         & spatial bound pili number density\\[-2.5pt] 
$\b{g}$                       & cell number density gradient                    &$\b{v}$                       & cell velocity\\[-2.5pt] 
$\b{y}$                       & non-linear deformation map                      &$\b{F}$                       & deformation gradient\\[-2.5pt]
$J$                           & Jacobian of the deformation gradient            &$\b{K}$                       & cofactor of the deformation gradient\\[-2.5pt] 
$\delta$                      & Kronecker delta                                 &$\Psi_\text{tot}$             & total internal energy \\[-2.5pt]
$\b{I}$                       & material second-order identity tensor           &$\b{i}$                       & spatial second-order identity tensor \\[-2.5pt] 
$\c{B}_0$                     & material configuration                          &$\p\c{B}_0$                   & boundary of the material configuration  \\[-2.5pt]
$\b{t}^{\text{a}}$            & active traction on material configuration       &$\b{t}^{\text{p}}$            & passive traction on material configuration  \\[-2.5pt]
$\delta \b{y}$                & linear momentum balance test function           &$\delta c$                    & cell number density conservation test function \\[-2.5pt]
$\delta \b{g}$                & cell density gradient continuity test function  &$\delta p$                    & bound pili number density evolution test function  \\[-2.5pt]
$\b{E}$                       & Green--Lagrange strain tensor                   &$\b{B}$                       & Piola deformation tensor \\[-2.5pt]
$\b{S}^{\text{a}}$            & active Piola--Kirchhoff stress                  &$\b{\tau}^{\text{a}}$         & active Kirchhoff stress \\[-2.5pt]
$\b{P}^{\text{a}}$            & active Piola stress                             &$\b{\sigma}^{\text{a}}$       & active Cauchy stress\\[-2.5pt]
$\b{P}^{\text{p}}$            & passive Piola stress                            &$\b{\sigma}^{\text{p}}$       & passive Cauchy stress\\[-2.5pt]
$\b{S}^{\text{f}}$            & pili formation Piola--Kirchhoff stress          &$\b{\sigma}^{\text{f}}$       & pili formation Cauchy stress\\[-2.5pt]
$\o{\b{P}}$                   & cell density gradient continuity Piola stresses &$\b{N}$                       & material unit normal to the boundary \\[-2.5pt]
$\v{R}$                       & assembled residual vector                       &$\b{U}$                       & global vector of unknowns  \\[-2.5pt]
$\v{K}$                       & assembled tangent stiffness                     &$\# e$                        & number of elements  \\[-2.5pt]
$N$                           & shape function of the finite elements           &$\v{I}^{\text{sym}}$          & symmetric fourth-order identity tensor \\[-2.5pt]
$\lambda$                     & penalty parameter                               &$\Delta t$                    & time step \\[-2.5pt]
\bottomrule
\end{tabular}
\end{table}

The remainder of this manuscript is organized as follows.
Table~\ref{tab:notation} gathers the key definitions and notations of the paper.
Section~\ref{sec:Governing equations} introduces the problem definition and presents the governing equations.
Finite element implementation of the problem is elaborated in Section~\ref{sec:FEM}.
Our proposed theory is illustrated through a set of numerical examples in Section~\ref{sec:result}.
Finally, Section~\ref{sec:conc} concludes the work and provides further outlooks.


\section{Governing equations}\label{sec:Governing equations}
\noindent
This section elaborates on the governing equations.
First the problem of cell aggregation is defined and all the parameters and their roles are introduced.
Afterwards, the continuum approach within the Lagrangian settings is detailed.
Note that for simplicity and readability, our formulations here are developed for a two-dimensional case, but there are no conceptual limitations to generalize it to three dimensions.

\subsection{Problem definition}\label{sec:prob}

\begin{figure}[b!]
\centering
{\includegraphics[width=0.6\textwidth]{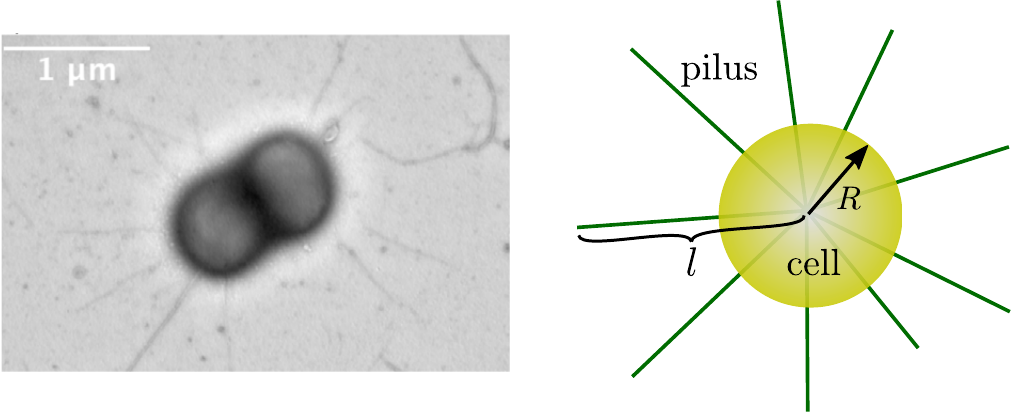}}
\caption{
Left: An image of N. gonorrhoeae bacterium and its pili obtained by transmission electron microscopy in~\cite{Kuan2021}.
Right: Simplification of the cell geometry to a circular shape for our analysis (right).
}
\label{fig:cell}
\end{figure}

\begin{figure}[t!]
\centering
{\includegraphics[width=1.0\textwidth]{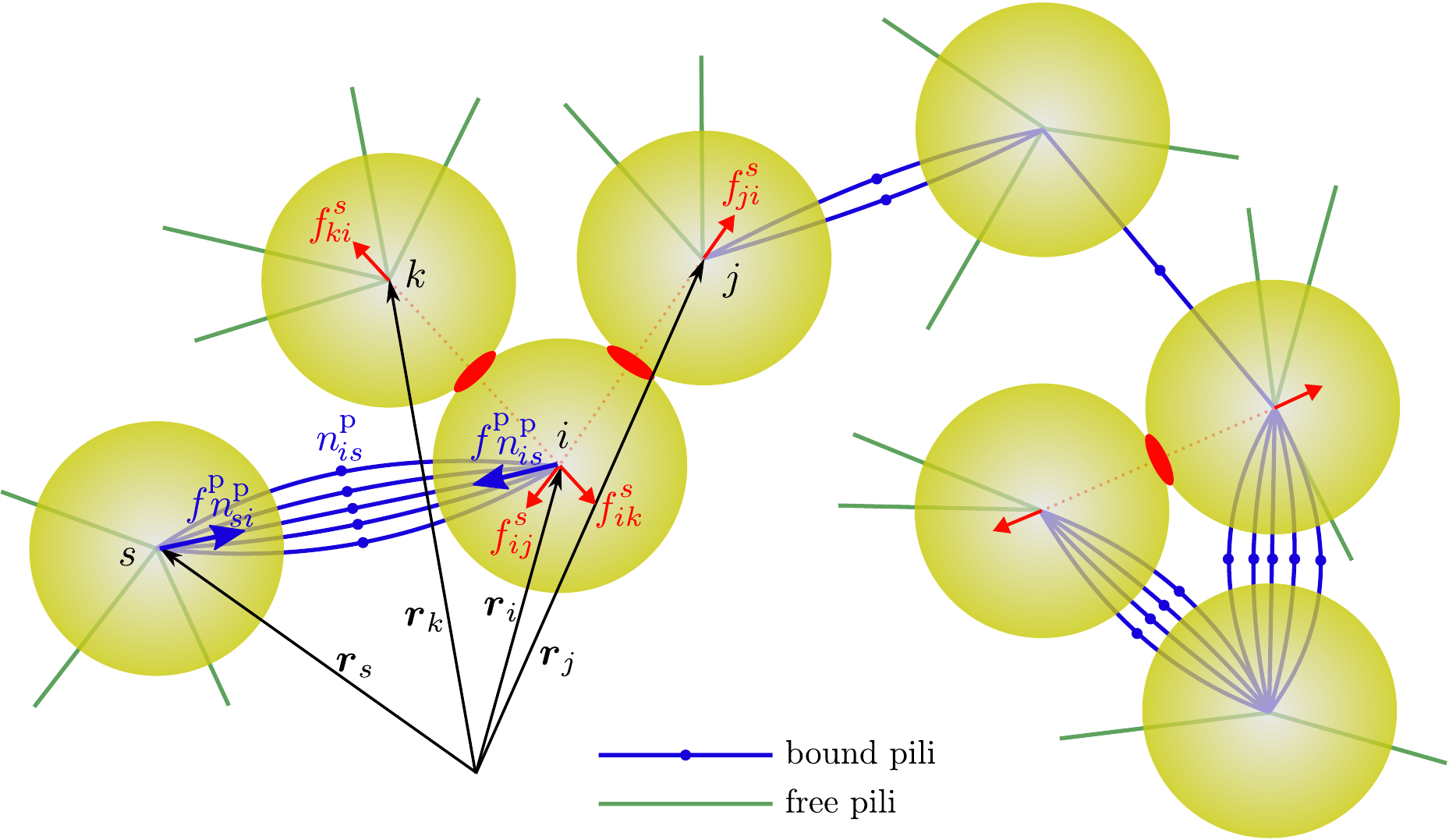}}
\caption{
A sketch of N. gonorrhoeae bacteria interacting by their pili.
The forces acting on each bacterium are illustrated.
Each bacterium $i$ is distinguished by its position vector $\b{r}_i$.
The inter-bacterium forces are steric repulsion forces between two attached bacteria which are shown in red, and the pili-pili mediated attractive forces between two bacteria that have formed a pili network are shown in blue.
These two forces act in the same direction of the line connecting the centers of two adjacent bacteria.
Bound pili are depicted by blue dotted lines whereas free pili are depicted by solid green lines.
Note that there could be multiple pili pairs pulling two neighboring cells together, as counted by $n_{is}^{\text{p}}$.
}
\label{fig:problem}
\end{figure}

\noindent
Figure~\ref{fig:cell} (left) shows a transmission electron microscopy image of a single NG bacterium together with its pili.
In this contribution, we approximate the cell geometry with a circular disk with the radius $R$ for the sake of simplicity, as shown in Fig.~\ref{fig:cell} (right).
Each NG bacterium is surrounded by approximately $10-20$ pili which are isotropically distributed around the cell~\cite{Zaburdaev2014,Eriksson2015}.
The length of each individual pilus $l$ was shown to be exponentially distributed with the average value of $\langle l \rangle = \ell_{0}$.
Measurements report the average length $\ell_{0}$ to be around $1-2\, \mu \text{m}$ with cell radius $R$ being around $1\,\mu \text{m}$~\cite{Holz2010}.

\begin{figure}[t!]
 \centering
 \includegraphics[width=1.0\textwidth]{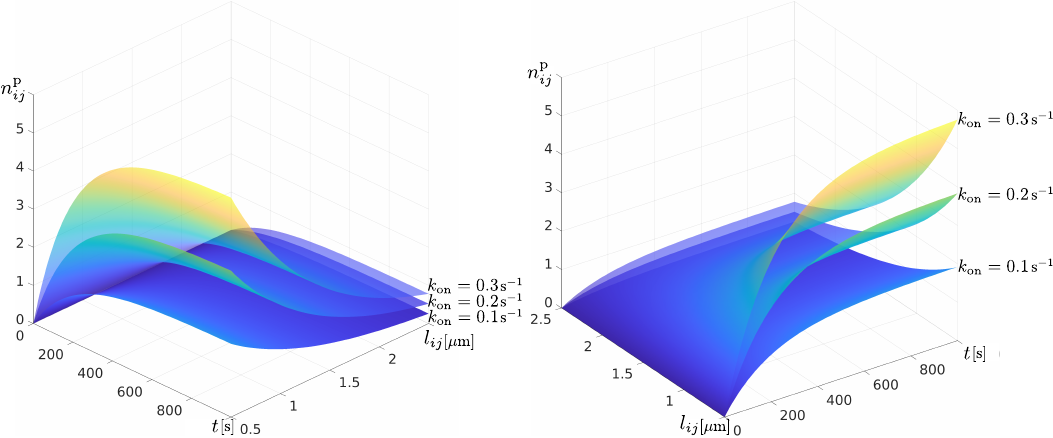}
 \caption{
 Variation of the bound pili pairs $n_{ij}^{\text{p}}$ between the cells $i$ and $j$ with respect the their distance $l_{ji}$ and time $t$ according to relation~\eqref{eq:pili_evolution}.
 Three different values $0.1$, $0.2$ and $0.3$ for the pili binding rate are considered.
 Two different views are provided for better illustration.
 Pili average length is set to $\ell_0=1.0\,\mu\text{m}$ and the unbinding rate is set to $k_{\text{off}}=0.005\,\text{s}^{-1}$.
 }
\label{fig:pili_evolution}
 \end{figure}

Figure~\ref{fig:problem} depicts a group of cells interacting via their pili.
In order to move on a substrate, cells use their pili.
Pili can grow, attach to the substrate and retract, which produces the force required for the displacement.
Pili-substrate interactions play a significant role in determining the cell motility which is essential for understanding the kinetics of aggregation.
Here, since we are mainly interested in examining the behavior of aggregates, only a substrate friction is considered to represent the pili-substrate interactions.
Further details on pili-substrate interactions are available in~\cite{Ponisch2019}.
The growth and retraction of the pili occurs through the process of polymerization and depolymerization which is powered by specific motor protein complexes~\cite{Marathe2014,Craig2019,Craig2004}.
Apart from attaching to the substrate, pili can also attach to pili of other cells.
Retraction of an attached pili pair generates an attractive force $f^\text{p}$ which pulls the cells towards each other.
In Fig.~\ref{fig:problem}, bound pili are depicted by blue lines with a dot on them, whereas free pili are shown by solid green lines.
The number of the pili pairs that have formed between two cells is denoted as $n^\text{p}$.
Accordingly, as illustrated for the cells $i$ and $s$, the overall attractive force between two cells which have formed a network of bound pili is $f^\text{p}n^\text{p}_{is}$.
Assuming no bound pairs at the initial time, the number of bound pili between the cells $i$ and $j$ (that do not move relative to each other) can be obtained via the relation 
\begin{equation}
\f{\d n_{ij}^{\text{p}}}{\d t} = \f{k_\text{on}\,e^{-l_{ij}\ell_0}}{2\pi\ell_0^2} - k_\text{off}n_{ij}^{\text{p}}
\qquad \Longrightarrow \qquad 
n_{ij}^{\text{p}} = \f{k_\text{on}\,e^{-l_{ij}\ell_0}}{2\pi\ell_0^2k_\text{off}}\left[ 1 - e^{-k_\text{off}\,t}   \right]\,,
\label{eq:pili_evolution}
\end{equation}
with $k_\text{on}$ being the pili binding rate and $k_\text{off}$ being the pili unbinding/detachment rate.
The distance between the cells $i$ and $j$ is $l_{ij}$ which is the magnitude of the vector pointing from cell $j$ to cell $i$ as $l_{ij}=|| \b{l}_{ij} ||=|| \b{r}_{i} - \b{r}_{j} ||$.
This relation was obtained by Kuan et al.~\cite{Kuan2021b,Kuan2021} which is a mean field approximation that ignores the discreteness of the pili number.
The factor $\f{k_\text{on}\,e^{-l_{ij}/\ell_0}}{2\pi\ell_0^2}$ stems from the assumption of an exponential distribution of the pili length and integration over all possible binding points along the line connecting two bacteria.
Figure~\ref{fig:pili_evolution} renders the variation of bound pili pairs between the cells $i$ and $j$ with respect to their distance $l_{ij}$ and time for three different values of binding rate $k_{\text{on}}$ (two different views are given for better illustration).
The average pili length is set to $\ell_0=1.0\,\mu\text{m}$.
It is confirmed that as time evolves, the number of bound pili pairs reaches a steady state.
Larger distance between the cells yields a decreased number of bound pili pairs.
This is justifiable since the further the cells are from each other, the less chance their pili have to attach to each other due to their limited length.
Evidently, larger binding rates also result in more bound pili pairs.

Another force that plays a role in the cell aggregate formation dynamics is the steric repulsive force $f^{\text{s}}_{ij}$ which occurs between the cells $i$ and $j$ that are in direct contact.
In the model, having this repulsive force is necessary to prevent interpenetration of the cells.
In summary, three major forces determine the dynamics of the cells in the network; cell-substrate friction, pili-pili mediated attractive force and steric repulsive force.
The forces arising from the contraction of the constantly remodeling pili network attract the cells together and in the presence of the excluded volume interactions, a densely packed colony is formed which behaves as an active visco-elastic material~\cite{Kuan2021}.

Since the pili-pili interactions are driven by active retractions and growth of pili, the system is inherently out of equilibrium.
The force balance equation for an individual cell in the aggregate can be obtained as
\begin{equation}
f^{\text{p}}\sum_{\substack{j=1\\j\neq i}} n_{ij}^{\text{p}}\hat{\b{l}}_{ji} + \sum_{\substack{j=1\\j\neq i}} f^{\text{s}}_{ij}\hat{\b{l}}_{ij} - \xi \b{v}_{i} = 0\,,
\label{eq:force-balance}
\end{equation}
with $\b{v}_i$ being the velocity of cell $i$ and $\xi$ being the friction coefficient between the cell and the substrate.
Note the first summation is over the neighboring cells with which the cell has formed pili pairs and the second summation is over the neighboring cells with which the cell is in direct contact.
The unit vector $\hat{\b{l}}_{ij}=\b{l}_{ij}/l_{ij}$ shows the direction of the pili-pili mediated force acting on cell $i$ by cell $j$ which is in the same direction of the steric repulsive force between these two cells.
The three terms in Eq.~\eqref{eq:force-balance} represent pili-pili mediated forces, excluded volume interactions and cell-substrate friction, respectively.
The summations in Eq.~\eqref{eq:force-balance} imply that the equation is solved for each cell individually.
While this approach offers more accuracy together with a higher degree of spatial resolution, it suffers from computational costs as the number of cells increases.
Moreover, if one is interested in examining the collective/overall behavior of the aggregate, tracking each individual cell is neither desirable nor necessary.
In this section, we propose a continuum model to describe the dynamics of dense bacterial colonies.
In our methodology, the system is characterized by two key features which distinguish it from previous models of active systems.
First, at time scales smaller than the pili detachment time, the bound pili network endows aggregates with elastic-like material properties for which a Lagrangian approach might be more suitable.
At larger time scales, pili can rearrange which allows for stress relaxation resulting in a fluid-like behavior for which an Eulerian approach is the common choice. 
Second, the attractive pili-pili mediated force dipoles are balanced by steric repulsion forces which allow for considering dense cellular aggregates.
These two key features are captured by the proposed continuum equations.


\subsection{Eulerian approach}\label{sec:Eulerian}
\noindent
In the Eulerian approach, the cell aggregate is treated as a compressible fluid-like material.
Therefore, our domain would be a fixed control window and we observe the changes that the aggregate goes through over this window.
This section briefly details on the governing equations of the problem within an Eulerian approach that were previously derived by coarse-graining microscopic equations in~\cite{Kuan2021}.

\subsubsection{Balance equations}\label{sec:balance_Eulerian}
\noindent
The first balance equation is the overall cell number conservation equation which reads
\begin{equation}
\begin{aligned}
\f{\p c_t}{\p t} +\nabla_{\!\b{x}}\!\cdot(c_t \b{v}) = 0\,,
\end{aligned}
\label{eq:continuity_Eulerian}
\end{equation}
where $c_t$ is the spatial cell \emph{number} density and $\b{v}$ is the cell (ether) velocity.
Note, the subscript $\b{x}$ indicates that spatial derivatives are taken with respect to the coordinates in the current configuration.
The two variables $c_t$ and $\b{v}$ are the unknowns of our problem.
The second balance equation is the linear momentum balance which arises from the microscopic balance of forces in Eq.~\eqref{eq:force-balance} and reads
\begin{equation}
\nabla_{\!\b{x}}\!\cdot\b{\sigma}^{\text{a}} + \nabla_{\!\b{x}}\!\cdot\b{\sigma}^{\text{p}} - \xi c_t\b{v} = \bnull\,,
\label{eq:LMB_Eulerian}
\end{equation}
with $\b{\sigma}^\text{a}$ and $\b{\sigma}^\text{p}$ being the active and passive Cauchy stresses, respectively.
The third balance equation is the bound pili number density evolution equation which reads
\begin{equation}
\f{\p p_{t}}{\p t}  + \nabla_{\!\b{x}}\! \cdot \lp p_{t}\b{v} \rp - k_{\text{on}}  \lb c_{t}^{2} + \f{3\ell_{0}^{2}}{4} c_{t} \,  \lb\Delta_{\b{x}} c_{t}\rb - \f{3\ell_{0}^{2}}{4} \big| \nabla_{\!\b{x}} c_{t} \big|^{2}  \rb + k_{\text{off}} p_{t} = 0\,,
\label{eq:pili_Eulerian}
\end{equation}
with $p_t$ the spatial bound pili number density.
This equation incorporates the pili turnover dynamics from Eq.~\eqref{eq:pili_evolution}.
Note, in our continuum formulation we denote the spatial pili number density as $p_t$ which is a coarse-grained version of $n_{ji}^{\text{p}}$.
The spatial Laplacian operator $\Delta_{\b{x}}$ is defined by the inner product $\Delta_{\b{x}}\sth\!:=\nabla_{\b{x}}\cdot\lp\nabla_{\b{x}}\sth\rp$.
Equation~\eqref{eq:pili_Eulerian} clearly states that the overall bound pili number is not conserved since the pili network constantly changes due to the binding and unbinding pairs.

\subsubsection{Stress definitions}\label{sec:stress_Eulerian}
\begin{figure}[h!]
\centering
\includegraphics[width=0.8\textwidth]{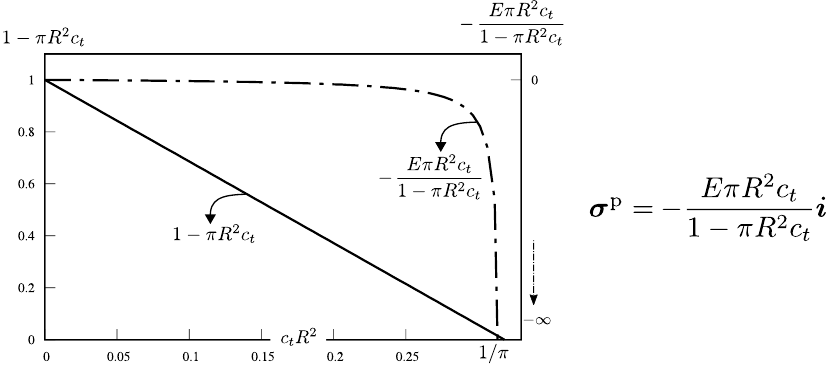}
\caption{
Illustration of the behavior of passive stress with respect to cell number density.
The denominator of the stress becomes zero at the value $\text{max}(c_t R^2)=1/\pi$ hence, the asymptotic behavior.
This value represents the upper bound for the cell number density.
}
\label{fig:piR}
\end{figure}

\noindent
The next step is to state the definitions for the active and passive stresses.
For the passive stress, we have a straightforward explicit expression which reads
\begin{equation}
\b{\sigma}^\text{p}=-\f{E \pi R^2 c_t}{1-\pi R^2 c_t}\b{i}\,,
\label{eq:passive_stress_Eulerian}
\end{equation}
with $\b{i}$ being the spatial second-order identity tensor, $E$ being the bulk modulus of the cell and $R$ being the cell radius.
The expression for the passive Cauchy stress is similar to a \emph{hydrostatic pressure-like stress}.
This expression is one of the simplest forms to describe the excluded volume interactions.
It resembles the van der Waals gas law for pressure in the absence of attractive interactions.
As a result, in the absence of the active stress, the only forces acting on cells would be steric repulsive forces and cell-substrate friction, thus while moving on the substrate, the cells tend to repel each other until they reach a state of equilibrium.
As it will be elucidated in the numerical examples, this gives rise to formation of a uniform distribution of cells throughout the domain.
Figure~\ref{fig:piR} depicts the behavior of the passive stress with respect to the cell number density.
The $x$-axis represents the dimensionless parameter $c_t R^2$ with the left $y$-axis representing $1-\pi R^2 c_t$ and the right $y$-axis representing the pressure term in the passive stress.
At $c_tR^2=1/\pi $, the stress renders an asymptotic behavior since the denominator tends to zero.
As a result, this value represents the upper bound for the cell number density.

The definition of the active Cauchy stress is given as coarse-grained average of the microscopic expression (for details see~\cite{Kuan2021}). 
For us, its time evolution equation is of interest which reads
\begin{equation}
\begin{aligned}
&\f{\p \b{\sigma}^\text{a}}{\p t}  + \nabla_{\!\b{x}}\!\cdot(\b{\sigma}^\text{a}\dyad\b{v}) - 2\left[ \b{l} \cdot \b{\sigma}^\text{a} \right]^{\text{sym}} = - \f{1}{\ell_0 p_t f^{\text{p}}} \left[\b{\sigma}^\text{a}\dyad\b{\sigma}^\text{a}\right]:\b{l}^{\text{sym}}  + \b{\sigma}^\text{f} - k_\text{off}\,\b{\sigma}^\text{a}\,,\\
\end{aligned}
\label{eq:active_stress_Eulerian}
\end{equation}
with $\b{l}$ being the spatial velocity gradient $\b{l}=\grad\b{v}$ and $\b{\sigma}^{\text{f}}$ being the pili formation Cauchy stress tensor which will be defined below.
The active Cauchy stress results from the pili-pili mediated attractive forces between the cells which pull the cells together and lead to the formation of an aggregate.
For the pili formation Cauchy stress an explicit formulation is achievable.
Due to the intrinsic nature of pili-pili interactions as force dipoles, this stress tensor has a nematic symmetry and reads
\begin{equation}
\b{\sigma}^{\text{f}} = \f{f^{\text{p}}k_\text{on}\ell_{0}}{2} \lb c_{t}^{2}\b{i} - \f{3\ell_{0}^{2}}{4} \lb  2\lb \nabla_{\b{x}} c_{t} \dyad \nabla_{\b{x}} c_{t} \rb  + \big| \nabla_{\b{x}} c_{t} \big|^{2} \b{i} \rb + \f{3\ell_{0}^{2}}{4} \lb c_{t} \lb\Delta_{\b{x}} c_{t}\rb\b{i} + 2c_{t}\nabla^{2}_{\b{x}}c_{t} \rb \rb\,,
\label{eq:pili_formation_stress_Eulerian}
\end{equation}
where $\nabla^{2}_{\b{x}}$ denotes the spatial Hessian operator defined by $\nabla^{2}_{\b{x}}\sth := \nabla_{\b{x}}\lp \nabla_{\b{x}}\sth\rp$.


\subsection{Lagrangian approach}\label{sec:Lagrangian}
\noindent
This section elaborates the pertinent equations within a Lagrangian framework, hence, all the equations are written in the material configuration.
In the Lagrangian approach, we follow the trajectory of the cells as time elapses.
Key geometric notions from nonlinear continuum kinematics are thereby the integration of the cell velocity into the cell nonlinear deformation map $\b{x}=\b{y}(\b{X},t)$ with $\b{x}$ and $\b{X}$ denoting the Eulerian and the Lagrangian coordinates.
Based on the notion of the cell deformation map $\b{y}$, we introduce the corresponding material deformation gradient $\b{F}:=\nabla_{\b{X}}\b{y}$, its Jacobian (determinant) $J:=\det \b{F}$ and its cofactor $\b{K}:=J\b{F}^{-T}$.
As will be elucidated, the main advantage of the Lagrangian approach is that it yields a significant simplification for the equations which facilitates finite element implementation of the problem.
In addition, employing the Lagrangian approach enables utilization of implicit time integration schemes which is computationally more robust compared to explicit time integration and less prone to instability issues.
A significant step towards the Lagrangian formulation of the problem is parametrization of all fields in Lagrangian coordinates $\b{X}$.
Especially the spatial cell number density is parameterized as $c_{t}=c_{t}(\b{X},t)$.
This critical step allows us to adopt $\c{C}^{0}$ continuous or discontinuous schemes and to avoid complications regarding the implementation of $\c{C}^{1}$ continuous elements.

\subsubsection{Balance equations}\label{sec:balance_Lagrangian}
\noindent
Similar to the previous section, we start with the overall cell number conservation equation.
Using the material time derivative $D\sth/Dt =\dot{\sth} = \p\sth/\p t + \nabla_{\b{x}}\sth \cdot \b{v}$ and Eq.~\eqref{eq:continuity_Eulerian}, the Lagrangian form of the overall cell number conservation equation can be written as
\begin{equation}
\f{D\lp Jc_{t}\rp}{D t}
=
\dot{\o{Jc_{t}}}
=
\dot{J} c_{t} + \dot{c_{t}} J
=0\,.
\label{eq:continuity_Lagrangian}
\end{equation}
The second balance equation is the linear momentum balance.
To define the Lagrangian version of the linear momentum balance equation, the velocity field is replaced by the time derivative of the deformation map as $\b{v}=\dot{\b{y}}$.
Using the relation $\nabla_{\b{X}}\cdot \b{P}=J \nabla_{\b{x}} \cdot \b{\sigma}$ the material linear momentum balance reads
\begin{equation}
\nabla_{\b{X}}\cdot\b{P}^\text{a} + \nabla_{\b{X}}\cdot\b{P}^\text{p} -\xi J c_{t}\dot{\b{y}} = \bnull\,,
\label{LMB_Lagrangian}
\end{equation}
with $\b{P}^\text{a}$ and $\b{P}^\text{p}$ being the active and passive Piola stresses, respectively.
The next step is to define the Lagrangian version of the bound pili evolution equation.
Using the relations $p_{0} =J p_{t}$ and $\dot{\sth} = \p\sth/\p t + \nabla_{\b{x}}\sth \cdot \b{v}$, the bound pili evolution equation can be written in terms of the material pili density as
\begin{equation}
\dot{p_{0}} - J k_\text{on}  \lb c_{t}^{2}  + \f{3\ell_{0}^{2}}{4} c_{t} \Delta_{\b{x}} c_{t} - \f{3\ell_{0}^{2}}{4} \big| \nabla_{\b{x}} c_{t} \big|^{2} \rb + k_{\text{off}}p_{0}=0\,.
\label{eq:pili_Lagrangian}
\end{equation}
Note, in Eq.~\eqref{eq:pili_Lagrangian}, the derivatives of the cell number density are still with respect to the spatial coordinates.
In Section~\ref{sec:gradient_enhanced_Lagrangian}, these derivatives will be transformed completely to the material configuration hence, our fully Lagrangian formulation.

\subsubsection{Stress definitions}\label{sec:stress_Lagrangian}
\noindent
Similar to the passive Cauchy stress, derivation of the passive Piola stress is straightforward.
Using the relation $\b{P} = J \b{\sigma}\cdot \b{F}^{-T} $, the passive Piola stress reads
\begin{equation}
\b{P}^\text{p}=-\f{E \pi R^2 c_{t}}{1-\pi R^2 c_{t}}\b{K}\,.
\label{passive_stress_Lagrangian}
\end{equation}
Derivation of the active Piola stress from the relation~\eqref{eq:active_stress_Eulerian} is an intricate task and requires further attention.
The active Piola stress tensor, a two-point tensor, must be converted to a fully material stress tensor.
In doing so, firstly the active Cauchy stress $\b{\sigma}^\text{a}$ is transformed to the active Kirchhoff stress $\b{\tau}^\text{a}$ as
\begin{equation}
\b{\tau}^\text{a} = J \b{\sigma}^\text{a}\,,
\end{equation}
with its material time derivative
\begin{equation}
\dot{\b{\tau}^\text{a}}
= 
J \left[ \f{\p \b{\sigma}^\text{a}}{\p t} + \nabla_{\b{x}} \!\cdot (\b{\sigma}^\text{a}\dyad\b{v}) \right]\,.
\end{equation}
Via multiplying Eq.~\eqref{eq:active_stress_Eulerian} by the Jacobian $J$, we arrive at a corresponding relation for Eq.~\eqref{eq:active_stress_Eulerian} in terms of the Kirchhoff stress
\begin{equation}
\dot{\b{\tau}^\text{a}} - 2\left[ \b{l} \cdot \b{\tau}^\text{a} \right]^{\text{sym}} = - \f{1}{\ell_0p_0f^{\text{p}}} \left[\b{\tau}^\text{a}\dyad\b{\tau}^\text{a}\right]:\b{l}^{\text{sym}}  + \b{\tau}^\text{f} - k_\text{off}\,\b{\tau}^\text{a}\,.
\label{eq:active_kirchhof_time_derivative}
\end{equation}
To obtain a fully material form of this equation, the Lie time derivative of the Kirchhoff stress needs to be introduced.
To express the Lie time derivative of the Kirchhoff stress, we need to calculate the pull-back of the Kirchhoff stress, calculate its material time derivative and then push it forward to the spatial configuration again which renders
\begin{equation}
\begin{aligned}
\c{L}_t\b{\tau}^\text{a} 
= 
\dot{\b{\tau}^\text{a}} - 2\left[ \b{l} \cdot \b{\tau}^\text{a} \right]^{\text{sym}}\,.
\end{aligned}
\label{eq:Lie_time_derivative}
\end{equation}
Note, the symmetry of the Kirchhoff stress is utilized in the above derivation.
Further details regarding the formulation of the Lie time derivative are available in~\ref{sec:app_2}.
The pull-back of the active Kirchhoff stress yields the active Piola--Kirchhoff stress as
\begin{equation}
\b{S}^\text{a} = \b{F}^{-1} \cdot \b{\tau}^\text{a} \cdot \b{F}^{-T}\,.
\end{equation}
Using Eqs.~\eqref{eq:active_kirchhof_time_derivative} and~\eqref{eq:Lie_time_derivative}, one can write
\begin{equation}
\c{L}_t\b{\tau}^\text{a} = \b{F} \cdot \dot{\b{S}^\text{a}}\cdot\b{F}^T = - \f{1}{\ell_0 p_0f^{\text{p}}} \left[\b{\tau}^\text{a}\dyad\b{\tau}^\text{a}\right]:\b{l}^{\text{sym}}  + \b{\tau}^\text{f} - k_\text{off}\,\b{\tau}^\text{a}\,.
\end{equation}
Accordingly, the material time derivative of the active Piola--Kirchhoff stress in terms of the active Kirchhoff stress is obtained as
\begin{equation}
\dot{\b{S}^\text{a}} 
= 
\b{F}^{-1} \cdot \lb\c{L}_t\b{\tau}^\text{a}\rb \cdot \b{F}^{-T}
=
\b{F}^{-1}\cdot\left[- \f{1}{\ell_0\rho_t^\text{p}f^{\text{p}}} \left[ \b{\tau}^\text{a}:\b{l}^{^{\text{sym}}}\right]\b{\tau}^\text{a}  + \b{\tau}^\text{f} - k_\text{off}\,\b{\tau}^\text{a}\right]\cdot\b{F}^{-T}\,.
\end{equation}
Finally, via replacing $\b{\tau}^\text{a}$ with $\b{F}\cdot\b{S}^\text{a}\cdot\b{F}^{T}$, we arrive at the final form of the material time derivative of the active Piola--Kirchhoff stress
\begin{equation}
\dot{\b{S}^\text{a}} = - \f{1}{\ell_0 p_0f^{\text{p}}} \left[ \b{S}^\text{a}:\dot{\b{E}}\right]\b{S}^\text{a}  + \b{S}^\text{f} - k_\text{off}\,\b{S}^\text{a}\,,
\label{eq:active_stress_Lagrangian}
\end{equation}
with $\b{E}$ being the Green--Lagrange strain tensor $\b{E}=\tfrac{1}{2}[\b{F}^T\cdot\b{F}-\b{I}]$.
For further details regarding the derivation of Eq.~\eqref{eq:active_stress_Lagrangian}, see~\ref{sec:app_2}.
The most significant advantage of using the Lagrangian approach is that it yields considerable simplification to the governing equations and the active stress time evolution equation.
The active stress can be numerically determined via time-discretizing Eq.~\eqref{eq:active_stress_Lagrangian}, see~\ref{sec:app_3} for further details.
Subsequently, after solving for the active Piola--Kirchhoff stress, the active Piola stress can be immediately obtained via $\b{P}^\text{a}=\b{F}\cdot\b{S}^\text{a}$.
Stating the material time derivative of the active Piola--Kirchhoff stress~\eqref{eq:active_stress_Lagrangian}, severely alleviates incrementally objective time integration of the active stress evolution.
The next step is to define the pili formation second Piola--Kirchhoff stress $\b{S}^{\text{f}}$.
Inserting the relation $\b{S} = J\b{F}^{-1}\cdot\b{\sigma}\cdot\b{F}^{-T}$ into Eq.~\eqref{eq:pili_formation_stress_Eulerian}, yields the pili formation Piola--Kirchhoff stress as
\begin{equation}
\begin{aligned}
\b{S}^\text{f} = \f{1}{2}f^{\text{p}}k_\text{on}\ell_{0}J \lb c_{t}^{2}\b{B} + \f{3\ell_{0}^{2}}{4} \lb  -2 \b{F}^{-1}\!\! \cdot \lb \nabla_{\b{x}} c_{t} \dyad \nabla_{\b{x}} c_{t} \rb \cdot \b{F}^{-T} \!\!\! - \big| \nabla_{\b{x}} c_{t} \big|^{2} \b{B} + c_{t}\lb \Delta_{\b{x}} c_{t}  \rb \b{B} + c_{t} \lb \b{F}^{-1}\!\! \cdot \nabla^{2}_{\b{x}}c_{t} \cdot \b{B} + \b{B} \cdot \nabla^{2}_{\b{x}}c_{t} \cdot \b{F}^{-T} \rb \rb \rb\,,\\
\end{aligned}
\label{eq:pili_formation_stress_Lagrangian}
\end{equation}
with $\b{B}=\b{F}^{-T}\cdot\b{F}^{-1}$ the Piola deformation tensor.
Note, in Eq.~\eqref{eq:pili_formation_stress_Lagrangian}, the gradients and divergences are with respect to the spatial coordinates which will be transformed completely to the material configuration
in Section~\ref{sec:gradient_enhanced_Lagrangian}.

\subsubsection{Weak form}\label{sec:weak_form_Lagrangian}
\noindent
To obtain the weak form of the overall cell number conservation equation, its strong form is multiplied by a scalar-valued test function $\delta c_{t}$ and then integrated over the referential domain as follows
\begin{equation}
\int_{\c{B}_0} \lb \dot{J} c_{t} + \dot{c_{t}} J \rb \delta c_{t} \, \d V 
=
0\, \qquad  \forall \delta c_{t}\,.
\label{eq:weak_continuity_Lagrangian}
\end{equation}
To obtain the weak form of the linear momentum balance, its strong form too is multiplied by a vector-valued test function $\delta \b{y}$ and then integrated over the referential domain which yields
\begin{equation}
 \int_{\c{B}_0}  \b{P}^{\text{a}} : \nabla_{\b{X}} \delta \b{y}  \, \d V + \int_{\c{B}_0}  \b{P}^{\text{p}} : \nabla_{\b{X}} \delta \b{y}  \, \d V  + \int_{\c{B}_0} \xi J c_{t} \dot{\b{y}} \cdot \delta \b{y} \, \d V = \int_{\p\c{B}_0}   \delta \b{y} \cdot \b{t}^{\text{a}}\, \d A + \int_{\p\c{B}_0}   \delta \b{y} \cdot \b{t}^{\text{p}}\, \d A\,\qquad  \forall \delta \b{y},
\label{eq:weak_LMB_Lagrangian_2}
\end{equation}
where $\b{t}^{\text{a}}$ and $\b{t}^{\text{p}}$ are active and passive tractions acting on the boundary of the domain, respectively.
Similarly, the weak form of the bound pili evolution equation is obtained via multiplying Eq.~\eqref{eq:pili_Lagrangian} by a scalar-valued test function $\delta p_0$ and then integrating over the referential domain as
\begin{equation}
\int_{\c{B}_0} \lb \dot{p_{0}} - J k_\text{on}  \lb c_{t}^{2}  + \f{3\ell_{0}^{2}}{4} c_{t} \Delta_{\b{x}} c_{t} - \f{3\ell_{0}^{2}}{4} \big| \nabla_{\b{x}} c_{t} \big|^{2} \rb + k_{\text{off}}p_{0} \rb \delta p_0 \,\d V= 0\,
\qquad  \forall \delta p_{0}.
\label{eq:weak_pili_Lagrangian}
\end{equation}

\subsubsection{Gradient enhanced framework}\label{sec:gradient_enhanced_Lagrangian}
\noindent
So far, the unknowns of our problem in Lagrangian approach have been the cell number density $c_{t}$, the cell deformation map $\b{y}$ and bound pili number density $p_{0}$.
A further challenge that arises in our problem is the presence of the Hessian of the spatial cell number density $\nabla^{2}_{\b{x}}c_{t}$ and the Laplacian of the spatial cell number density $\Delta_{\b{x}}c_{t}$ in the bound pili number density evolution equation~\eqref{eq:pili_Lagrangian} and the pili density formation Piola stress equation~\eqref{eq:pili_formation_stress_Lagrangian}.
Thus, second derivatives of the spatial cell number density $c_t$ are required to be calculated.
In doing so, different strategies have been employed in the literature, among which the well-established ones are $\c{C}^{1}$ continuous elements~\cite{Papanicolopulos2012,Fischer2012}, isogeometric analysis~\cite{Hughes2005,Bazilevs2010}, micromorphic continuum approach and the gradient enhanced framework~\cite{Kirchner2005,Neff2014,Kaessmair2020}.
In order to stay within the realm of the classical finite element method (FEM) associated with a $\c{C}^0$ continuous interpolation approach, we adopt the gradient enhanced framework.
In doing so, we introduce an additional independent spatial vector field $\b{g}$ to represent the spatial gradient of the cell number density $\nabla_{\b{x}} c_t$ and \emph{weakly enforce} $\b{g}\overset{!}{=}\nabla_{\b{x}} c_t$.
Note, similar to the spatial cell number density, its gradient is also parametrized as a function of the Lagrangian coordinates $\b{g}=\b{g}(\b{X},t)$.
The first step is to transform the Laplacian and Hessian to substitute the terms including spatial derivatives of $c_{t}$ with $\b{g}$ as
\begin{equation}
\begin{aligned}
&\Delta_{\b{x}} c_t = \nabla_{\b{x}} \cdot (\nabla_{\b{x}} c_t) = \nabla_{\b{x}} \cdot \b{g} \,,\\
&\nabla^{2}_{\b{x}} c_t = \nabla_{\b{x}} (\nabla_{\b{x}} c_t) = \nabla_{\b{x}}^{\text{sym}} \b{g} = \f{1}{2} \lb \nabla_{\b{x}} \b{g} + \nabla_{\b{x}}^{T} \b{g} \rb\,,
\end{aligned}
\label{eq:cell_gradient_Eulerian}
\end{equation}
where a symmetric gradient of $\b{g}$ follows from the symmetry of the Hessian of the cell number density.
Afterwards the resulting terms must be pulled back to the material configuration so as to unify the Lagrangian formalism.
For an arbitrary vector field $\b{a}$, the following relations hold between the spatial and material divergence and gradient $\nabla_{\b{x}} \cdot \b{a}  = \nabla_{\b{X}} \b{a} : \b{F}^{-T}$ and $\nabla_{\b{x}} \b{a} = \nabla_{\b{X}} \b{a} \cdot \b{F}^{-1}$. 
Accordingly, the Laplacian and Hessian of the cell number density can be written in terms of $\b{g}$ in the material configuration as
\begin{equation} 
\begin{aligned}
&\Delta_{\b{x}} c_{t} =  \nabla_{\b{x}} \cdot \b{g} = \nabla_{\b{X}} \b{g} : \b{F}^{-T}\,,\\
&\nabla^{2}_{\b{x}} c_{t} = \nabla_{\b{x}}^\text{sym} \b{g} = \f{1}{2} \lb \nabla_{\b{X}} \b{g} \cdot \b{F}^{-1}  + \b{F}^{-T} \cdot \nabla_{\b{X}}^{T} \b{g}  \rb\,.
\end{aligned}
\label{eq:cell_gradient_Lagrangian}
\end{equation}
Using Eq.~\eqref{eq:cell_gradient_Lagrangian}, the strong~\eqref{eq:pili_Lagrangian} and weak~\eqref{eq:weak_pili_Lagrangian} forms of the bound pili number density evolution equation can be stated as fully Lagrangian equations as
\begin{equation}
\dot{p_{0}} - J k_\text{on}  \lb c_{t}^{2}  + \f{3\ell_{0}^{2}}{4} c_{t} \lb \nabla_{\b{X}} \b{g}:\b{F}^{-T} \rb - \f{3\ell_{0}^{2}}{4} \big| \b{g} \big|^{2} \rb + k_{\text{off}}p_{0}=0\,,
\label{eq:pili_Lagrangian_updated_1}
\end{equation}
and 
\begin{equation}
\int_{\c{B}_0} \lb \dot {p_{0}} - J k_\text{on}  \lb c_{t}^{2}  + \f{3\ell_{0}^{2}}{4} c_{t} \lb \nabla_{\b{X}} \b{g}:\b{F}^{-T} \rb - \f{3\ell_{0}^{2}}{4} \big| \b{g} \big|^{2} \rb + k_{\text{off}}p_{0} \rb \delta p_0 \, \d V=0\,,
\label{eq:pili_Lagrangian_updated_2}
\end{equation}
respectively.
Accordingly the pili density formation Piola--Kirchhoff stress equation~\eqref{eq:pili_formation_stress_Lagrangian} can be rewritten as
\begin{equation}
\begin{aligned}
\b{S}^\text{f} = \f{1}{2}f^{\text{p}}k_\text{on}\ell_{0}J \lb c_{t}^{2}\b{B} + \f{3\ell_{0}^{2}}{4} \lb  -2 \b{F}^{-1}\!\! \cdot \lb \b{g} \dyad \b{g} \rb \cdot \b{F}^{-T} \!\!\! - \lb \b{g} \rb^{2} \b{B} + c_{t}\lb \nabla_{\b{X}} \b{g}:\b{F}^{-T}  \rb\b{B} + c_{t} \lb \b{F}^{-1}\!\! \cdot \nabla_{\b{X}}\b{g} \cdot \b{B} + \b{B} \cdot \nabla_{\b{X}}^{T}\b{g} \cdot \b{F}^{-T} \rb \rb \rb\,.\\
\end{aligned}
\label{eq:pili_formation_stress_Lagrangian_updated}
\end{equation}

\subsubsection{Fictitious variational approach}\label{sec:variational_Lagrangian}
\noindent
The final step to complete our formulation is to update the governing equations given the new continuity requirement.
The advantage of using a (fictitious) variational approach is that it provides the admissible forms of tractions and external forces without any prior knowledge or assumptions.
In addition, with no extra effort or derivations, the variational approach furnishes the weak form which is central to the finite element implementation of the problem.
In order to obtain the updated governing equations, a total energy functional is minimized via setting its first variation to zero.
To begin, we assume that there exist a fictitious energy function $\Psi_1$ whose variation with respect to the degrees of freedom yields the already existing weak forms~\eqref{eq:weak_continuity_Lagrangian}, \eqref{eq:weak_LMB_Lagrangian_2} and \eqref{eq:pili_Lagrangian_updated_2}.
The term ``fictitious'' for this energy function implies that the original weak forms were directly derived from the strong forms and not from any specific energy function and, as mentioned before, the strong forms themselves are obtained using the coarse-grained continuum approach developed in~\cite{Kuan2021,Kuan2021b}.
To enforce the continuity of the spatial cell number density gradient, we introduce a new energy function $\Psi_2$ whose definition will be provided shortly.
The total energy $\Psi_{\text{tot}}$ consists of the initial fictitious energy and the newly defined energy
\begin{equation}
\Psi_{\text{tot}} = \Psi_1 + \Psi_2\,.
\end{equation}
The energies are the integrals of their corresponding internal energy densities over their associated domains as
\begin{equation}
\Psi_{\text{tot}} = \int_{\c{B}_0}  \psi_{\text{tot}} \, \d V \,, \qquad \Psi_1 = \int_{\c{B}_0} \psi_1 \, \d V\,, \qquad \Psi_2 = \int_{\c{B}_0} \psi_2 \, \d V \,.
\end{equation}
To minimize $\Psi_{\text{tot}}$, its first variation is set to zero.
That is
\begin{equation}
\delta \Psi_{\text{tot}} \stackrel{\b{.}}{=} 0 \qquad  \Longrightarrow \qquad  \delta\Psi_1 + \delta\Psi_2 \stackrel{\b{.}}{=} 0\,.
\end{equation}
The field variables in our problem are the cell number density $c_t$, cell deformation map $\b{y}$, bound pili number density $p_0$ and cell number density gradient $\b{g}$.
The energy density function to impose the continuity of the spatial cell density gradient reads
\begin{equation}
\psi_2 = \psi_2\lp c_t,\b{y},\b{g} \rp = \f{1}{2} J \lambda \lb \b{F}^{-T}\cdot \nabla_{\b{X}}c_{t} - \b{g} \rb^2 \,,
\label{eq:energy_Lagrangian}
\end{equation}
with $\lambda$ being the penalty parameter.
The penalty parameter determines how strongly the condition $\b{g}\overset{!}{=}\nabla_{\b{x}} c_t$ is satisfied.
Note, the dependence of $\psi_2$ on $\b{y}$ is through the deformation gradient $\b{F}$ and its determinant $J$.
Accordingly, one could write
\begin{equation}
\int_{\c{B}_0}\!\! \delta \psi_{\text{tot}} \, \d V 
=
\int_{\c{B}_0} \delta \psi_1 \, \d V +  \int_{\c{B}_0} \delta \lp \f{1}{2} J \lambda \lb \b{F}^{-T}\cdot \nabla_{\b{X}}c_{t} - \b{g} \rb^2 \rp \, \d V\,.
\label{eq:variation_Lagrangian_1}
\end{equation}
Calculating the variation of the overall energy density with respect to all degrees of freedom and their gradients, yields the weak forms as
\begin{equation} 
\int_{\c{B}_0} \delta \psi_{\text{tot}} \, \d V
=
\int_{\c{B}_0}\!\!\! \delta \psi_1 \, \d V 
+
\int_{\c{B}_0}\!\!\! \delta \psi_2 \, \d V\,,
\end{equation}
with
\begin{equation}
\begin{aligned}
&\int_{\c{B}_0} \delta \psi_{1} \, \d V
=
\int_{\c{B}_0}\!\!\! \delta_{c_{t}} \psi_1 \, \d V
+
\int_{\c{B}_0}\!\!\! \delta_{\b{y}} \psi_1 \, \d V
+
\int_{\c{B}_0}\!\!\! \delta_{p_0} \psi_1 \, \d V \,,\\
&\int_{\c{B}_0} \delta \psi_{2} \, \d V
=
\int_{\c{B}_0}\!\!\! \delta_{c_{t}} \psi_2 \, \d V
+
\int_{\c{B}_0}\!\!\! \delta_{\b{y}} \psi_2 \, \d V
+
\int_{\c{B}_0}\!\!\! \delta_{\b{g}} \psi_2 \, \d V \,.
\end{aligned}
\end{equation}
The variations of the newly defined energy function read
\begin{equation} 
\begin{aligned}
&\int_{\c{B}_0}\! \delta_{c_{t}} \psi_2  \, \d V = \int_{\c{B}_0} \lambda J \lb \b{B} \cdot \nabla_{\b{X}} c_{t} - \b{F}^{-1} \cdot \b{g} \rb \cdot  \nabla_{\b{X}} \delta c_{t} \, \d V\,, \\[5pt]
&\int_{\c{B}_0}\! \delta_{\b{y}} \psi_2   \, \d V = \int_{\c{B}_0}\! \underbrace{\lb\f{1}{2}\lambda J\lb \b{F}^{-T} \cdot \nabla_{\b{X}} c_{t} - \b{g} \rb^{2} \b{F}^{-T} - \lambda J \lb \lb \b{F}^{-T}\cdot \nabla_{\b{X}} c_{t}  \rb \dyad \lb \b{B}\cdot \nabla_{\b{X}} c_{t}  - \b{F}^{-1} \cdot \b{g}\rb\rb \rb }_{\o{\b{P}}} : \nabla_{\b{X}} \delta\b{y} \, \d V \,,\\[5pt] 
&\int_{\c{B}_0}\! \delta_{\b{g}} \psi_2  \, \d V = - \int_{\c{B}_0} \lambda J \lb \b{F}^{-T} \cdot \nabla_{\b{X}} c_{t} - \b{g} \rb \cdot \delta \b{g} \, \d V \,,
\end{aligned}
\end{equation}
with $\o{\b{P}}$ being the cell-density-gradient-continuity-induced Piola stress which arises due to the dependence of $\psi_2$ on the deformation.
Finally, we arrive at the final weak form of the governing equations which form our residual system as
\begin{equation}
\begin{aligned}
&\text{overall cell number conservation:} &&\int_{\c{B}_0} \!\!\! \lb \dot{J} c_{t} + \dot{c_{t}} J \rb \delta c_{t} \, \d V + \int_{\c{B}_0} \!\!\! \lambda J \lb \b{B} \cdot \nabla_{\b{X}} c_{t} - \b{F}^{-1}\cdot\b{g} \rb \cdot  \nabla_{\b{X}} \delta c_{t} \, \d V = 0\,,\\
&\text{linear momentum balance:} &&\int_{\c{B}_0} \!\!\! \xi J c_{t} \dot{\b{y}} \cdot \delta \b{y} \, \d V +\!\!\! \int_{\c{B}_0} \!\!\! \lb \b{P}^{\text{a}} + \b{P}^{\text{p}} + \o{\b{P}} \rb :  \nabla_{\b{X}} \delta\b{y} \, \d V  - \int_{\p\c{B}_0} \!\!\!\!  \delta \b{y} \cdot \lb \b{t}^{\text{a}} + \b{t}^{\text{p}} \rb\, \d A=\bnull  \,, \\
&\text{bound pili number density evolution:} && \int_{\c{B}_0} \lb \dot {p_{0}} - J k_\text{on}  \lb c_{t}^{2}  + \f{3\ell_{0}^{2}}{4} c_{t} \lb \nabla_{\b{X}} \b{g}:\b{F}^{-T} \rb - \f{3\ell_{0}^{2}}{4} \big| \b{g} \big|^{2} \rb + k_{\text{off}}p_{0} \rb \delta p_0 \, \d V = 0\,, \\
&\text{cell number density gradient continuity:} &&\int_{\c{B}_0} \!\!\! -\lambda J \lb \b{F}^{-T}\cdot \nabla_{\b{X}} c_{t} - \b{g} \rb \cdot \delta \b{g} \, \d V = \bnull\,.
\end{aligned}
\label{eq:weak_Lagrangian}
\end{equation}

\section{Finite element implementation}\label{sec:FEM}
\noindent
In this section, we present a general finite element formulation for the implementation of our developed methodology.
The first step towards finite element implementation of our problem is the derivation of the discretized weak forms which has been concluded in Section~\ref{sec:variational_Lagrangian}.
For time integration, the time interval $\mathds{T}$ is subdivided into a set of intervals $\Delta t$ with
\begin{equation}
\mathds{T} = \bigcup_{n=0}^{\#\text{ts}-1} \lb t_n,\, t_{n+1}  \rb\,,
\end{equation}
where $\#\text{ts}$ denotes the number of time steps and the time increment is defined by $\Delta t = t_{n+1}-t_{n}$.
As mentioned earlier, an implicit time integration associated with the backward Euler method is employed here in which an unknown $y_{n+1}$ is determined via the relation $y_{n+1}-y_{n} =  f(y_{n+1}) \Delta t$ based on the continuous evolution equation $\dot{y}=f(y)$.

Here, the spatial discretization of the problem domain is carried out using the Bubnov--Galerkin finite element method.
The geometry together with the variables are approximated using isoparametric coordinates $\b{\xi} \in [-1, 1]^2$.
Using standard interpolations together with the isoparametric concept, the geometry is discretized as
\begin{equation}
\begin{aligned}
\b{X}\Big|_{\c{B}_0}  \!\!\!\approx \b{X}(\b{\xi}) = \sum_{i=1}^{\#\text{e}} N^{i}(\b{\xi}) \b{X}^i\,, 
\end{aligned}
\label{eq:discretized_Lagrangian}
\end{equation}
with $\# \text{e}$ denoting the number of elements and $\b{\xi}$ denoting the natural space coordinates.
Accordingly the discretized fields read
\begin{equation}
\begin{aligned}
&c_{t}\Big|_{\c{B}_0} \!\!\!\approx c_{t}(\b{\xi}) = \sum_{i=1}^{\#\text{e}} N^{i}(\b{\xi}) c_{t}^i\,, 
&&\delta c_{t}\Big|_{\c{B}_0} \!\!\!\approx \delta c_{t}(\b{\xi}) = \!\!\!\sum_{i=1}^{\#\text{e}} N^{i}(\b{\xi}) \delta c_{t}^i\,,\\
&\b{y}\Big|_{\c{B}_0} \!\!\!\approx \b{y}(\b{\xi}) = \sum_{i=1}^{\#\text{e}} M^{i}(\b{\xi}) \b{y}^i\,, 
&&\delta\b{y}\Big|_{\c{B}_0}  \!\!\!\approx \delta\b{y}(\b{\xi})  = \!\!\!\sum_{i=1}^{\#\text{e}} M^{i}(\b{\xi}) \delta\b{y}^i\,, \\
&p_{0}\Big|_{\c{B}_0} \!\!\!\approx p_{0}(\b{\xi}) = \sum_{i=1}^{\#\text{e}} N^{i}(\b{\xi}) p_{0}^i\,, 
&&\delta p_{0}\Big|_{\c{B}_0} \!\!\!\approx \delta p_{0}(\b{\xi}) = \!\!\!\sum_{i=1}^{\#\text{e}} N^{i}(\b{\xi}) \delta p_{0}^i\,,\\ 
&\b{g}\Big|_{\c{B}_0} \!\!\!\approx \b{g}(\b{\xi}) = \sum_{i=1}^{\#\text{e}} N^{i}(\b{\xi}) \b{g}^i\,,
&&\delta\b{g}\Big|_{\c{B}_0}  \!\!\!\approx \delta\b{g}(\b{\xi})  = \!\!\!\sum_{i=1}^{\#\text{e}} N^{i}(\b{\xi}) \delta\b{g}^i\,,
\end{aligned}
\label{eq:discretized_Lagrangian}
\end{equation}
with $N^i$ and $M^i$ being different shape functions, which, however, can also be chosen identical in particular cases.
For the cases where a mixed finite element formulation is required, these shape functions are polynomials of different orders.

The discretized weak form of the overall cell number conservation equation~$\eqref{eq:weak_Lagrangian}_1$ reads 
\begin{equation}
\int_{\c{B}_0} \lb c_{t_{n}} \f{ J_{n+1} - J_{n} }{\Delta t} + J_{n+1} \f{ c_{t_{n+1}} - c_{t_{n}} }{\Delta t} \rb   N^{i}  \, \d V + \int_{\c{B}_0} \lambda J_{n+1} \lb \b{B}_{n+1} \cdot \nabla_{\b{X}} c_{t_{n+1}} - \b{F}_{n+1}^{-1} \cdot \b{g}_{n+1} \rb \cdot  \nabla_{\b{X}}  N^{i}  \, \d V
=0\,.
\label{eq:discretized_continuity_Lagrangian}
\end{equation}
The discretized weak form of the linear momentum balance equation~$\eqref{eq:weak_Lagrangian}_2$ reads
\begin{equation}
\int_{\c{B}_0} \xi J_{n+1} c_{t_{n+1}} \f{\b{y}_{n+1}-\b{y}_{n}}{\Delta t} M^{i}  \, \d V + \int_{\c{B}_0} \lb \b{P}^{\text{a}}_{n+1} + \b{P}^{\text{p}}_{n+1} + \o{\b{P}}_{n+1} \rb \cdot \nabla_{\b{X}}  M^{i} \, \d V 
=
\int_{\p\c{B}_0} \!\!\!  \b{t}^{\text{a}}M^{i}\, \d A +\!\!\! \int_{\p\c{B}_0} \!\!\!  \b{t}^{\text{p}}M^{i}\, \d A\,.
\label{eq:discretized_LMB_Lagrangian}
\end{equation}
The discretized weak form of the bound pili number density evolution equation~$\eqref{eq:weak_Lagrangian}_3$ reads 
\begin{equation}
\int_{\c{B}_0} \lb \f{p_{0_{n+1}}-p_{0_{n}}}{\Delta t} - J_{n+1} k_\text{on}  \lb c_{t_{n+1}}^{2}  + \f{3\ell_{0}^{2}}{4} c_{t_{n+1}} \lb \nabla_{\b{X}} \b{g}_{n+1}:\b{F}_{n+1}^{-T} \rb - \f{3\ell_{0}^{2}}{4} \big| \b{g}_{n+1} \big|^{2} \rb + k_{\text{off}}p_{0_{n+1}} \rb N^{i} \, \d V = 0\,,
\label{eq:discretized_pili_Lagrangian}
\end{equation}
and finally, the discretized weak form of the cell density gradient continuity equation~$\eqref{eq:weak_Lagrangian}_4$ reads
\begin{equation}
\int_{\c{B}_0} -\lambda J_{n+1} \lb \b{F}_{n+1}^{-T} \cdot \nabla_{\b{X}} c_{t_{n+1}} - \b{g}_{n+1} \rb  N^i \, \d V
=
\bnull\,.
\label{eq:discretized_gradient_Lagrangian}
\end{equation}
The fully-discrete form of balance equations can be obtained by applying the spatial approximation which form a residual system associated with the global node $I$ as
\begin{equation}
\begin{aligned}
&\v{R}^{I}_{c} = \assembly{\alpha=1}{\#e} \int_{\c{B}_0} \lb c_{t_{n}} \f{ J_{n+1} - J_{n}}{\Delta t} + J_{n+1} \f{ c_{t_{n+1}} - c_{t_{n}}}{\Delta t}  \rb   N^{i}  \, \d V + \assembly{\alpha=1}{\#e} \int_{\c{B}_0} \lambda J_{n+1} \lb \b{B}_{n+1} \cdot \nabla_{\b{X}} c_{t_{n+1}} - \b{F}_{n+1}^{-1}\cdot\b{g}_{n+1} \rb \cdot  \nabla_{\b{X}}  N^{i}  \, \d V = 0\,,\\
&\v{R}^{I}_{\b{y}} = \assembly{\alpha=1}{\#e} \int_{\c{B}_0} \xi J_{n+1} c_{t_{n+1}} \f{\b{y}_{n+1}-\b{y}_{n}}{\Delta t} M^{i}  \, \d V + \assembly{\alpha=1}{\#e}\int_{\c{B}_0} \lb \b{P}^{\text{a}}_{n+1} + \b{P}^{\text{p}}_{n+1} + \o{\b{P}}_{n+1} \rb \cdot \nabla_{\b{X}}  M^{i} \, \d V - 
\assembly{\alpha=1}{\#se}\int_{\p\c{B}_0} \!\!\! \lb \b{t}^{\text{a}} + \b{t}^{\text{p}} \rb M^{i}\, \d A= \bnull \,,\\
&\v{R}^{I}_{p} = \assembly{\alpha=1}{\#e} \int_{\c{B}_0} \lb \f{p_{0_{n+1}}-p_{0_{n}}}{\Delta t} - J_{n+1} k_\text{on}  \lb c_{t_{n+1}}^{2}  + \f{3\ell_{0}^{2}}{4} c_{t_{n+1}} \lb \nabla_{\b{X}} \b{g}_{n+1}:\b{F}_{n+1}^{-T} \rb - \f{3\ell_{0}^{2}}{4} \big| \b{g}_{n+1} \big|^{2} \rb + k_{\text{off}}p_{0_{n+1}} \rb N^{i} \, \d V = 0\,,\\
&\v{R}^{I}_{\b{g}} = \assembly{\alpha=1}{\#e} \int_{\c{B}_0} -\lambda J_{n+1} \lb \b{F}_{n+1}^{-T} \nabla_{\b{X}} c_{t_{n+1}} - \b{g}_{n+1} \rb  N^i \, \d V = \bnull \,,
\end{aligned}
\label{eq:residual_Lagrangian}
\end{equation}
with $\# se$ denoting the number of surface elements and $\assembly{}{}$ being the assembly operator.
The global residual vector consist of the above four residual vectors as
\begin{equation}
\v{R}_{\text{tot}}=
\lb
\begin{matrix}
\v{R}_{c_t}\\
\v{R}_{\b{y}}\\
\v{R}_{p_0}\\
\v{R}_{\b{g}}
\end{matrix}
\rb\,,
\qquad 
\text{with} 
\qquad
\v{R}_{c_t} = 
\lb
\begin{matrix}
\v{R}^{1}_{c_t}\\
\v{R}^{2}_{c_t}\\
\vdots\\
\v{R}^{\text{nn}}_{c_t}\\
\end{matrix}
\rb\,,
\qquad
\v{R}_{\b{y}} = 
\lb
\begin{matrix}
\v{R}^{1}_{\b{y}}\\
\v{R}^{2}_{\b{y}}\\
\vdots\\
\v{R}^{\text{nn}}_{\b{y}}\\
\end{matrix}
\rb\,,
\qquad
\v{R}_{p_0} = 
\lb
\begin{matrix}
\v{R}^{1}_{p_0}\\
\v{R}^{2}_{p_0}\\
\vdots\\
\v{R}^{\text{nn}}_{p_0}\\
\end{matrix}
\rb\,,
\qquad
\v{R}_{\b{g}} = 
\lb
\begin{matrix}
\v{R}^{1}_{\b{g}}\\
\v{R}^{2}_{\b{g}}\\
\vdots\\
\v{R}^{\text{nn}}_{\b{g}}\\
\end{matrix}
\rb\,,
\label{eq:global_residual_Lagrangian}
\end{equation}
where nn denotes the total number of nodes.
Finally, the fully discrete nonlinear system of governing equations becomes
\begin{equation}
\v{R}_{\text{tot}} = \v{R}_{\text{tot}} (\b{U}) \stackrel{\b{!}}{=} \bnull\,,
\label{eq:residual_zero_Lagrangian}
\end{equation}
with $\b{U}$ being the global vector of unknowns.
To find the solution of the system~\eqref{eq:residual_zero_Lagrangian}, the Newton--Raphson scheme is employed.
The consistent linearization of the resulting system of equations yields
\begin{equation}
\v{R}_{\text{tot}} (\b{U}_{n+1}) = \v{R}_{\text{tot}} (\b{U}_{n}) + \v{K}_{\text{tot}}\cdot \Delta\b{U}_{n} \stackrel{\b{!}}{=} \bnull
\qquad
\text{with}
\qquad
\v{K}_{\text{tot}} = \f{\p \v{R}_{\text{tot}}}{\p \b{U}}\big\lvert_{n}
\qquad 
\text{and}
\qquad
\b{U}_{n+1}=\b{U}_{n}+\Delta\b{U}_{n}\,,
\label{eq:Newton_Lagrangain}
\end{equation}
where the subscript $n$ indicates the step number.



\section{Numerical results}\label{sec:result}
\noindent
This section aims to illustrate our proposed theory through a set of numerical examples.
In doing so, the evolution of the cell number density under various conditions is investigated.
Moreover, parametric studies are carried out in order to highlight the influence of different parameters on the evolution of the cell number density.
Throughout all the examples, the domain is a $80\times80$ square subject to periodic boundary conditions.
Implicit time integration is adopted with no stability condition which leads to robust solution of our problem and larger time steps could be used.
All the numerical results are obtained from our in-house finite element code.
The solution procedure is robust and for all examples, we obtain convergence with a quadratic rate associated with the Newton--Raphson scheme.


\subsection{Cell number density evolution}\label{sec:active}

\begin{figure}[b!]
\centering
\includegraphics[width=1.0\textwidth]{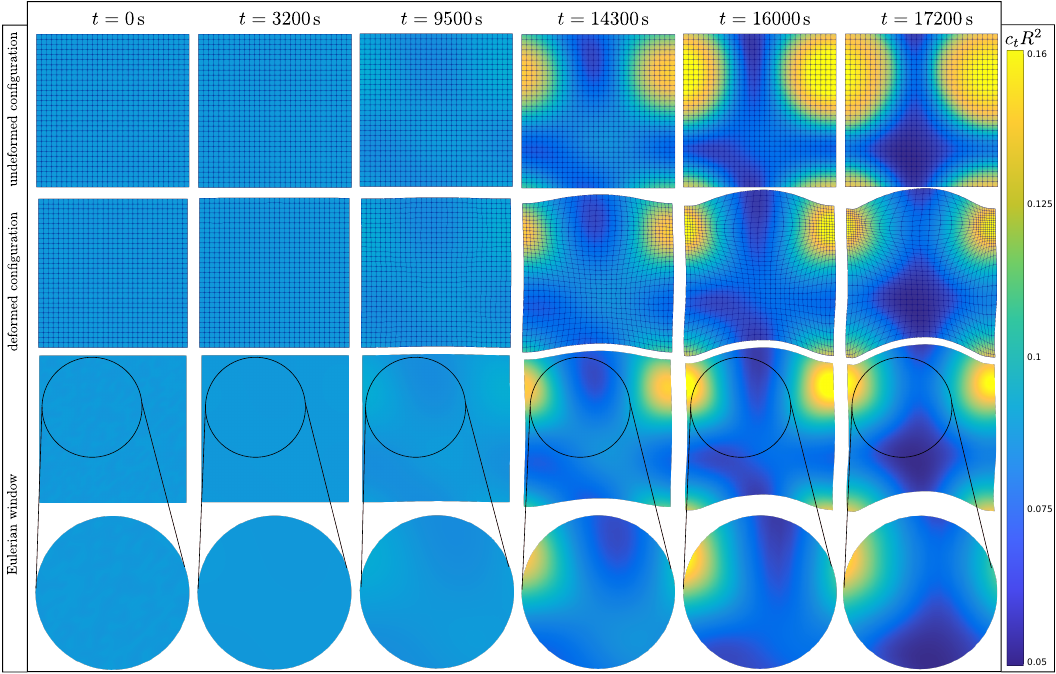}
\caption{
Snapshots of the cell number density evolution in the presence of the active stress.
The first and second rows represent the cell number density distribution in undeformed and deformed configurations, respectively.
The third and the fourth rows illustrate how an Eulerian observer would see the cell number density evolution through a fixed window associated with the Eulerian framework.
}
\label{fig:active}
\end{figure}

\begin{figure}[b!]
\centering
\includegraphics[width=1.0\textwidth]{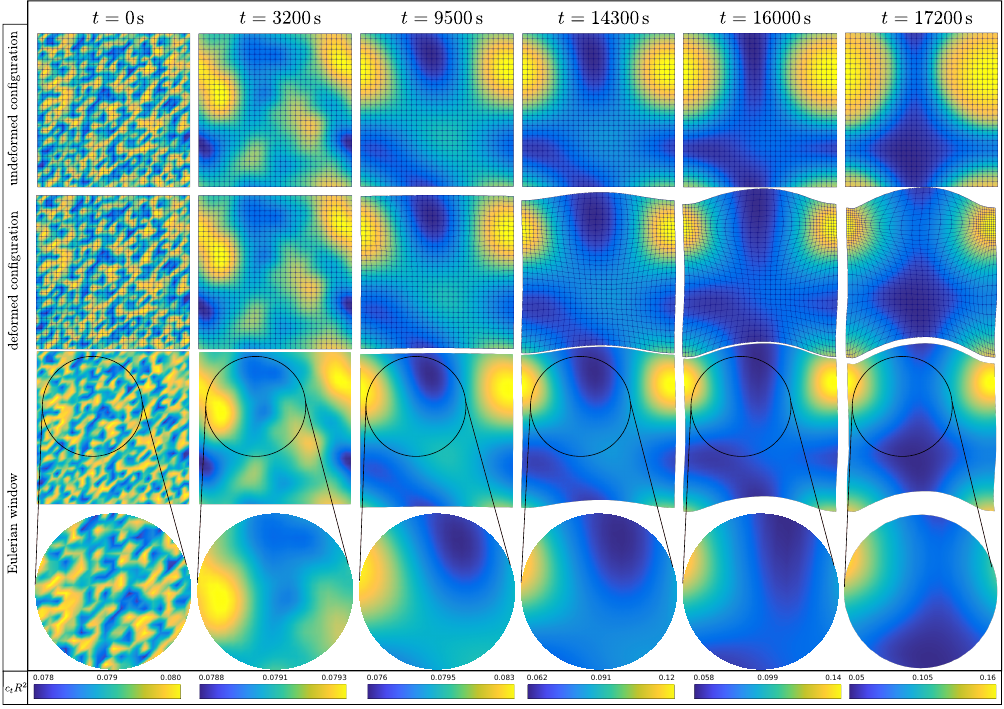}
\caption{
Another illustration of the cell aggregate evolution in Fig.~\ref{fig:active} with variable color-bar range.
}
\label{fig:active_2}
\end{figure}

\noindent
This section investigates the evolution of the cells under two different scenarios.
In the first scenario, the cellular aggregate behavior is investigated in the presence of all the driving forces.
That is, the pili-pili mediated attractive forces, steric repulsion forces and the cell-substrate friction.
In the second scenario, the problem is simplified to the case where the pili-pili mediated attractive forces are eliminated.
In this case, we investigate the cell number density evolution when only the steric repulsion forces and the cell-substrate friction are the forces acting on the cells.

Figure~\ref{fig:active} renders six different snapshots of the cellular aggregate evolution at different times associated with the first scenario.
The color represents the magnitude of cell number density throughout the domain.
The first and second rows represent the cell number density distribution in the undeformed and deformed configurations, respectively.
The third and the fourth rows illustrate how an Eulerian observer would see the cellular aggregate evolution through a fixed window associated with the Eulerian framework.
At the initial time $t=0$, the problem starts with a uniform distribution of cells with initial $c_t=0.079$ perturbed with $\pm 0.001$ relative random fluctuation.
The cell radius is set to $R=1\,\mu \text{m}$, the cell bulk modulus is set to $E=1\,\text{N}/\text{m}^2$ and the friction coefficient is set to $\xi=10\,\text{Ns}/\text{m}$.
In addition, the pili pair binding rate is $k_{\text{on}}=0.0178\,\text{s}^{-1}$, the pili pair unbinding rate is $k_{\text{off}}=0.01\,\text{s}^{-1}$, the average pili length is $\ell_0=2.0\,\mathrm{\mu} \text{m}$ and the pili-pili mediated attractive force is $f^{\text{p}}=18.0\,\text{pN}$.
Since the backward Euler time integration is adopted, an adaptive time stepping is employed with initial time step of $\Delta t=0.5\,\text{s}$.
If convergence is obtained at very few iterations, the time step is enlarged by a factor of $1.2$ and vice versa.
Through the aggregate evolution, first an initial homogeneous smooth cell distribution is observed due to the action of the passive stress.
This is justifiable since the formation of bound pili pairs (and thus emergence of active stress) requires time hence, the dominance of the passive stress at the starting steps.
As time elapses, more bound pili pairs form as depicted in Fig.~\ref{fig:pili_evolution}.
Thus cells tend to attract each other which is reflected in the dominance of the active stress which triggers the onset of the phase separation.
Therefore, pili-mediated attractive forces act as the major driving force which lead to the formation of an aggregate.
Figure~\ref{fig:active_2} provides another illustration of the cell aggregate evolution in Fig.~\ref{fig:active} with variable color range in each step.
The purpose of the figure is to elucidate the initial random distribution of cells and highlight the early smooth distribution of cells due to the passive stress which is reflected in the tighter color bar range.

\begin{figure}[t!]
\centering
\includegraphics[width=0.32\textwidth]{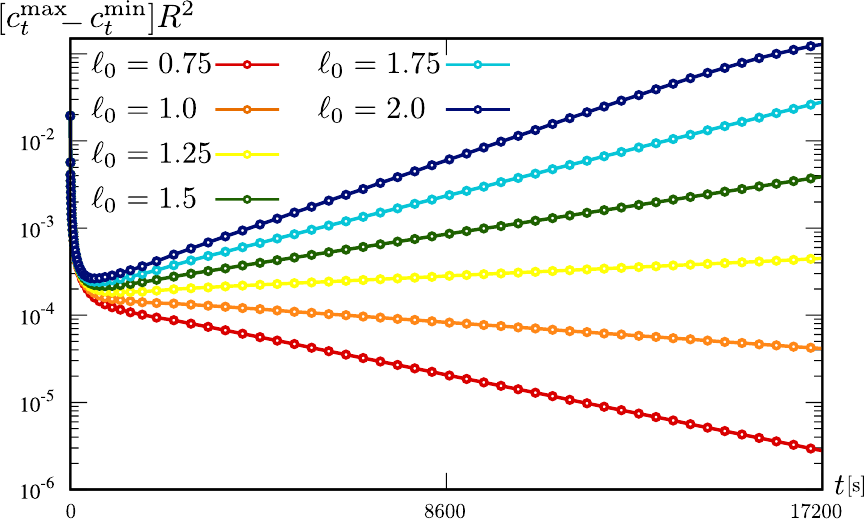}
\includegraphics[width=0.32\textwidth]{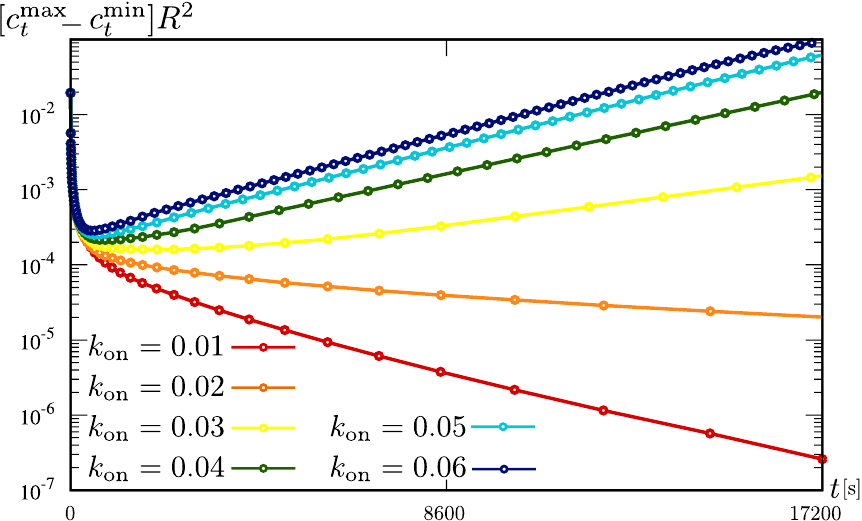}
\includegraphics[width=0.32\textwidth]{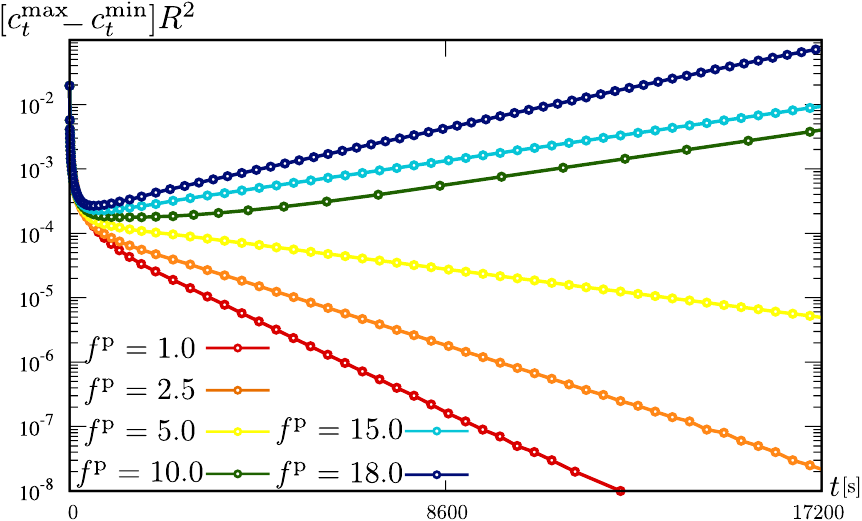}
\caption{
A parametric study on the cell density evolution.
The $y$-axis shows the difference between the maximum and minimum cell density throughout the domain while the $x$-axis shows the elapsed time.
The left figure studies the effects of the pili average length $\ell_0$ where $k_{\text{on}}=0.05\,\text{s}^{-1}$ and $f^{p}=12\,\text{pN}$.
The middle figure studies the effects of the pili pair binding rate $k_{\text{on}}$ where $\ell_0=2.0\,\mu \text{m}$ and $f^{p}=12\,\text{pN}$.
The right figure studies the effects of the pili pair attractive force $f^{p}$ where $k_{\text{on}}=0.05\,\text{s}^{-1}$ and $\ell_0=2.0\,\mu \text{m}$.
}
\label{fig:active_parameters}
\end{figure}

In Fig.~\ref{fig:active_parameters}, a parametric study is carried out to investigate the influence of the average pili length $\ell_0$, pili-pili binding rate $k_{\text{on}}$ and pili-pili mediated attractive force $f^{\text{p}}$ on the formation of an aggregate.
Each figure renders the difference between the maximum and the minimum cell number density versus time.
The left figure studies the effects of the pili average length $\ell_0$ where $k_{\text{on}}=0.05\,\text{s}^{-1}$ and $f^{p}=12\,\text{pN}$, the middle figure studies the effects of the pili pair binding rate $k_{\text{on}}$ where $\ell_0=2.0\,\mu \text{m}$ and $f^{p}=12\,\text{pN}$ and the right figure studies the effects of the pili pair attractive force $f^{p}$ where $k_{\text{on}}=0.05\,\text{s}^{-1}$ and $\ell_0=2.0\,\mu \text{m}$.
For all figures, an initial decrease in the cell number density difference is observed which is associated with function of the passive stress.
This behavior is vividly observed in the second step in Fig.~\ref{fig:active_2} as the color bar limits tend to tighten.
It is observed that, if the values of $f^{p}$, $\ell_0$ and $k_{\text{on}}$ are small, the pili-pili mediated forces cannot overcome the steric repulsive forces and the difference continues to decrease until a uniform homogeneous distribution is obtained.
For these cases the phase separation does not occur and the dynamics of the cell network is mainly driven by the steric repulsive forces and cell-substrate friction.
The lower the values of $f^{p}$, $\ell_0$ and $k_{\text{on}}$, the faster the equilibrium state is reached.
However, if the values of $f^{p}$, $\ell_0$ and $k_{\text{on}}$ are large enough, the pili-pili mediated forces tend to dominate the other forces as time elapses which is reflected in the larger cell number density difference.
In these cases, the active stress plays a more considerable role compared to the passive stress which yields the onset of the phase separation.
Increasing the values of $f^{p}$, $\ell_0$ and $k_{\text{on}}$ results in higher rates of phase separation as indicated by larger slopes of graphs in Fig.~\ref{fig:active_parameters}.
This behavior is natural to expect as larger values of $f^{p}$ means that the attractive force between the pili pairs is stronger which causes the cells to form a colony at a higher rate.
Larger values of $\ell_0$ indicate longer pili implying that pili are more capable to reach out and attach to pili of other cells which facilitates the formation of a colony.
And finally, larger values of $k_{\text{on}}$ signify faster binding rate between the pili leading to a faster formation of bound pairs which accelerates the colony formation.
Note, the parameters in the graphs in Fig~\ref{fig:active_parameters}, totally follow the condition obtained by Kuan et al.~\cite{Kuan2021,Kuan2021b} using linear stability analysis which reads
\begin{equation*}
-\f{1}{\xi}\f{E\pi R^2}{\lb1-c_t^{\text{initial}}\pi R^2 \rb^2}+\f{c_t^{\text{initial}}\ell_0f^{\text{p}}k_{\text{on}} }{k_{\text{off}}\xi}>0\,.
\end{equation*}
Satisfying the above condition indicates the onset of the phase separation.

\begin{figure}[b!]
\centering
\includegraphics[width=0.8\textwidth]{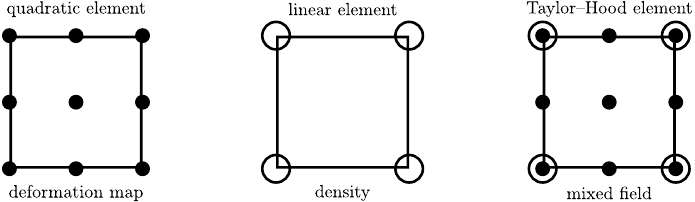}
\caption{Illustration of the Taylor--Hood element utilized for the simulation.
For this element type we used quadratic approximation for the deformation map field and linear approximation for the cell number density in order to observe the LBB condition.}
\label{fig:element}
\end{figure}

Now we consider the second scenario where the cell number density evolution is examined while its dynamics is driven only by the steric repulsive forces and the pili-substrate friction.
In the absence of the the pili-pili mediated attractive forces, the complexity of the problem reduces considerably since the bound pili number density evolution equation is eliminated from the system of equations.
In addition, there is no need to adopt a gradient enhanced framework to deal with the Laplacian and Hessian of the cell number density.
As a result, the field variables simply become the original cell number density $c_t$ and the deformation map $\b{y}$ and, under periodic boundary condition, the fully discrete residual system~\eqref{eq:residual_Lagrangian} simplifies to 
\begin{equation}
\begin{aligned}
&\v{R}^{I}_{c} = \assembly{\alpha=1}{\#e} \int_{\c{B}_0} \lb c_{t_{n}}  \f{J_{n+1} - J_{n}}{\Delta t} + J_{n+1}  \f{c_{t_{n+1}} - c_{t_{n}}}{\Delta t}  \rb   N^{i}  \, \d V = 0\,,\\
&\v{R}^{I}_{\b{y}} = \assembly{\alpha=1}{\#e} \int_{\c{B}_0} \xi J_{n+1} c_{t_{n+1}} \f{\b{y}_{n+1}-\b{y}_{n}}{\Delta t} M^{i}  \, \d V + \assembly{\alpha=1}{\#e}\int_{\c{B}_0} \b{P}^{\text{p}}_{n+1}  \cdot \nabla_{\b{X}}  M^{i} \, \d V  = \bnull\,.\\
\end{aligned}
\end{equation}
As mentioned earlier, pili-pili attractive forces are absent and the repulsive forces between the cells are reflected in the passive stress.
To cope with the mixed nature of the problem, i.e. to respect the LBB inf-sup condition, we allow interpolation of the cell number density $c_t$ and the deformation map $\b{y}$ with different polynomial order.
We exploit the Taylor--Hood element for our finite element implementation which is linear in $c_t$ and quadratic in $\b{y}$.
Subsequently, the corresponding shape functions $N^{i}$ and $M^{i}$ are polynomials of different orders.
Figure~\ref{fig:element} provides an illustrative schematic of the Taylor--Hood elements adopted in our simulations.

\begin{figure}[b!]
\centering
\includegraphics[width=1.0\textwidth]{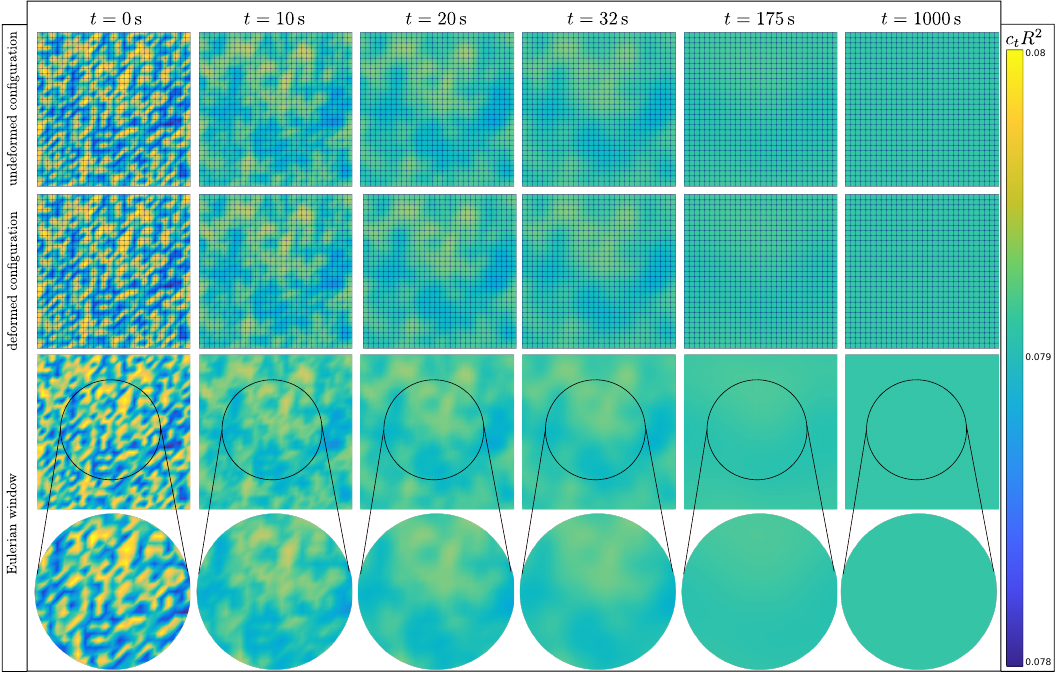}
\caption{
Different snapshots of the cell network evolution in the absence of the active stress.
The first and second rows represent the cell number density distribution in undeformed and deformed configurations, respectively.
The third and the fourth rows illustrate how an Eulerian observer would see the cell number density evolution through a fixed window associated with the Eulerian framework.
The dynamics of cell is only driven by the steric repulsion forces and the cell-substrate friction.
Uniform color throughout the domain indicates uniform homogeneous distribution of the cells.
}
\label{fig:passive}
\end{figure}

Figure~\ref{fig:passive} renders six different snapshots of the cell density evolution at different times.
The boundary and the initial conditions are similar to the previous case study.
Similarly the deformed and undeformed bodies are provided in the first and second rows where the third and fourth rows are associated with the Eulerian window.
Since there exists no pili-pili attractive force, as time elapses, the cells tend to repel each other until they reach to the state of equilibrium.
This equilibrium state is obtained when the cells are uniformly distributed throughout the domain.
It is observed that throughout time evolution, the difference between the cell densities tend to vanish which signifies the uniformity of the cell distribution.
Moreover, due to the absence of the active stress, the domain does not undergo any considerable deformation where the deformation gradient remains close to identity.

\begin{figure}[h!]
\centering
\includegraphics[width=0.48\textwidth]{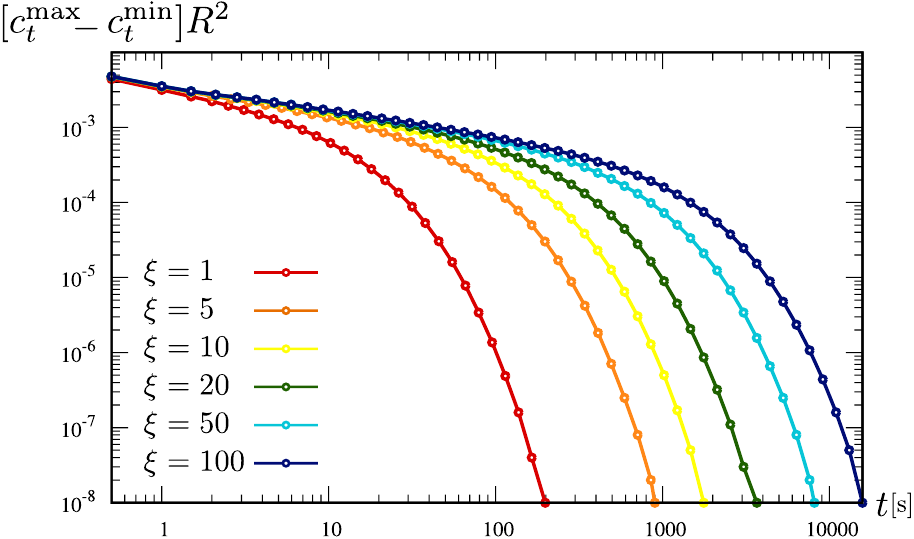}
\includegraphics[width=0.48\textwidth]{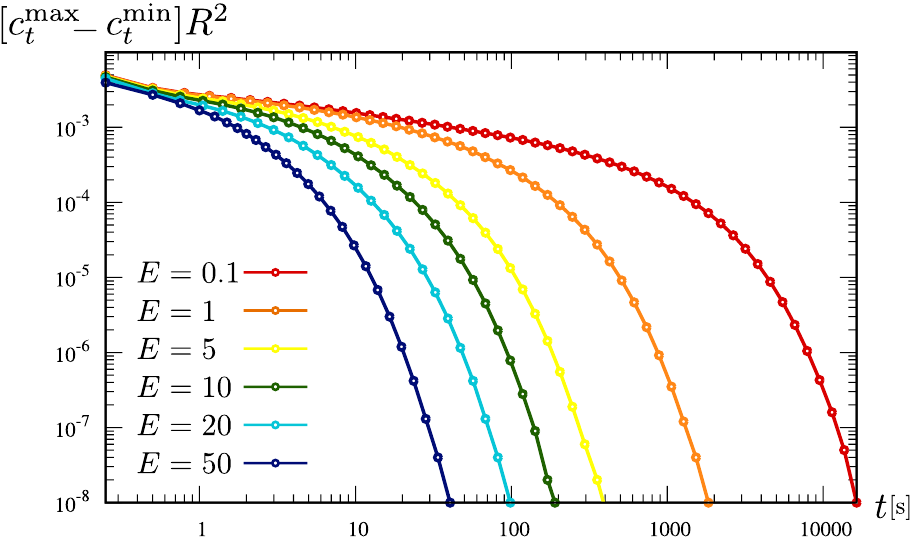}
\caption{
A parametric study on the cell density evolution in the absence of the active stress.
The $y$-axis shows the difference between the maximum and minimum cell density throughout the domain as a measure of homogeneity of the cell distribution while the $x$-axis shows the elapsed time.
The left figure renders the effects of friction coefficient on the cell density evolution when $E=1\,\text{N}/\text{m}^2$ and $R=1\,\mu \text{m}$.
The right figure renders the effects of cell bulk modulus on the cell density evolution when $R=1\,\mu \text{m}$ and $\xi=10\,\text{Ns}/\text{m}$.
}
\label{fig:passive_parameters}
\end{figure}

In Fig.~\ref{fig:passive_parameters}, a parametric study is carried out to investigate the influence of the friction coefficient $\xi$ and cell bulk modulus $E$ on the cell number density evolution.
Each figure renders the difference between the maximum and the minimum cell number density versus time.
The left figure studies the effects of the friction coefficient with $E=1\,\text{N}/\text{m}^2$ and $R=1\,\mu \text{m}$ and the right figure studies the effects of the cell bulk modulus with $\xi=10\,\text{Ns}/\text{m}$ and $R=1\,\mu \text{m}$.
Smaller density difference in each figure indicates more uniform distribution.
It is observed that increasing the friction coefficient delays the uniform distribution which is understandable since the friction coefficient impedes the cells movement.
Thus, for larger friction coefficients, the cells require more time to disperse due to the repulsive forces.
On the other hand, we observe that a uniform distribution is reached at a faster rate for larger cell bulk modulus.
This is also justifiable since larger bulk modulus indicates more rigidity of the cells which gives rise to their quicker dispersion after colliding with other cells.

\subsection{Overall cell number conservation}\label{sec:number-density} 
\begin{figure}[h!]
\centering
\includegraphics[width=1.0\textwidth]{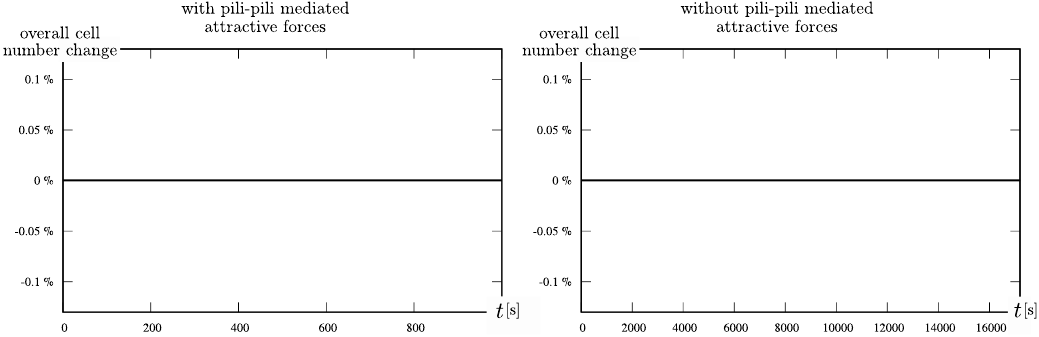}
\caption{
Illustration of the discrete conservation of the overall cell number for both case studies associated with Figs.~\ref{fig:active} and~\ref{fig:passive}.
The $x$ axis represents the elapsed time whereas the $y$ axis represents the percentage of change in the overall cell number.
The left figure corresponds to the case where all forces are present in the dynamics of the system whereas the left figure corresponds to the case where pili-pili mediated forces are absent.
}
\label{fig:mass}
\end{figure}

\noindent
It is noteworthy to mention that our mixed four-field $\lp c_t,\b{y},p_0,\b{g}\rp$ implementation conserves the total cell number in a Lagrangian solution domain exactly.
To ensure the discrete conservation of the total cell number throughout the domain, we rewrite the strong form of the continuity equation~\eqref{eq:continuity_Lagrangian}, integrate it over the referential domain and expand it as
\begin{equation}
\int_{\c{B}_0} \dot{\o{Jc_{t}}} \, \d V
=
\int_{\c{B}_0} \lb \dot{J} c_{t} + \dot{c_{t}} J \rb \, \d V
\approx
\int_{\c{B}_0} \lb \f{J_{n+1}-J_{n}}{\Delta t} c_{t_{n}} + \f{c_{t_{n+1}}-c_{t_{n}}}{\Delta t} J_{n+1} \rb \, \d V
=
0\,.
\label{eq:cell_density_conserve_1}
\end{equation}
Note, to satisfy the overall cell number conservation equation in the last integral, the cell density in the first term must be chosen from step $n$ whereas the Jacobian in the second term must be chosen from step $n+1$.
Accordingly, Eq.~\eqref{eq:cell_density_conserve_1} reads
\begin{equation}
\int_{\c{B}_0} J_{n+1}c_{t_{n+1}} \d V - \int_{\c{B}_0} J_{n}c_{t_{n}} \, \d V
=
0\,,
\label{eq:cell_density_conserve_2}
\end{equation}
which can in turn be written as the integrals over the spatial domain as
\begin{equation}
\int_{\c{B}_{t_{n+1}}} c_{t_{n+1}} \d v_{n+1}
=
\int_{\c{B}_{t_{n}}} c_{t_{n}} \d v_{n}
\,,
\label{eq:cell_density_conserve_3}
\end{equation}
which is equivalent to the conservation of the total cell number.
Figure~\ref{fig:mass} illustrates the percentage of the overall cell number change versus time for the two case studies associated with Figs.~\ref{fig:active} and~\ref{fig:passive}.
It is of great significance to point out that no cell number loss occurs within our systems, thus the variation of the total cell number remains zero without any fluctuation, see~\cite{Onate2014,Cremonesi2017} for further discussions on  the issue regarding mass loss/production in Lagrangian formulations.

\subsection{Colony coalescence}\label{sec:colony}

\begin{figure}[t!]
\centering
\includegraphics[width=1.0\textwidth]{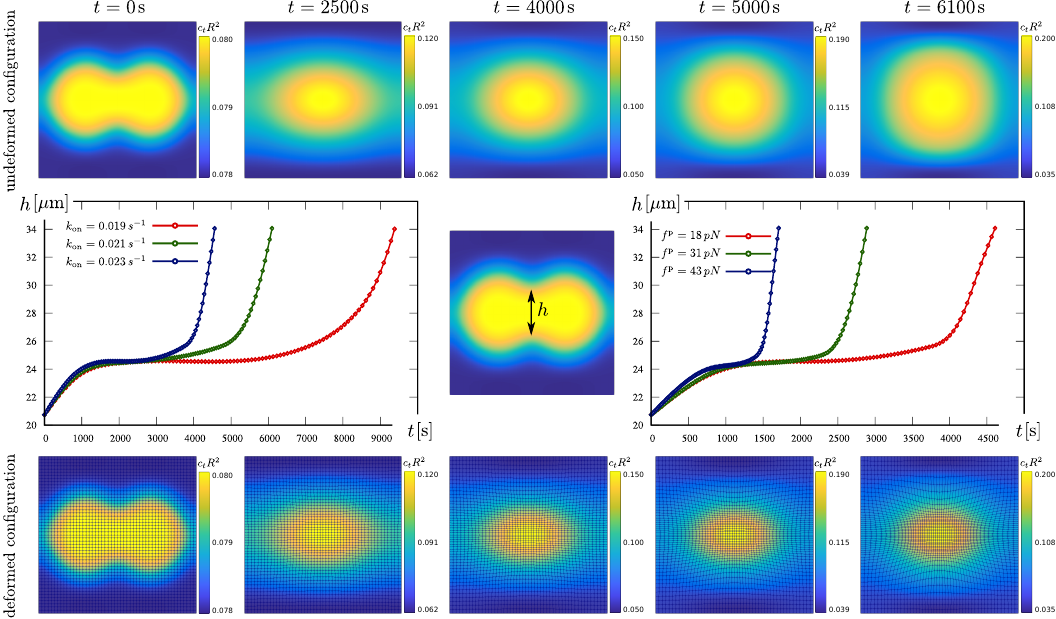}
\caption{
Illustration of the colony coalescence process.
Five snapshots of the coalescence process in both undeformed and deformed configurations are depicted.
The parameter $h$ measures the bridge length between the two merging colonies.
The center plots render the effects of $k_{\text{on}}$ and $f^{\text{p}}$ on the evolution of the bridge between the colonies.
}
\label{fig:two-bubble}
\end{figure}

\noindent
In this section, two colonies are put next to each other and the process of their coalescence is examined.
In doing so, we define a new parameter $h$ which measures the bridge length between the two merging colonies.
From a computational point of view, we define the boundary of the bridge where the the cell number density gradient in $y$ direction is $96\%$ of its maximum value throughout the domain.
Figure~\ref{fig:two-bubble} investigates the influence of the pili-pili binding rate $k_{\text{on}}$ and the pili-pili mediated attractive force $f^{\text{p}}$ on the size of the bridge.
Five snapshots of the coalescence process in both undeformed and deformed configurations are shown on top and bottom, respectively.
The snapshots correspond to the green line in the left plot with $k_{\text{on}}=0.021\,\text{s}^{-1}$ and $f^{\text{p}}=18\,\text{pN}$.
The left plot renders the bridge length versus time for three different values of the pili-pili binding rate with $f^{\text{p}}=18\,\text{pN}$.
The blue line corresponds to $k_{\text{on}}=0.019\,\text{s}^{-1}$, the green line corresponds to $k_{\text{on}}=0.021\,\text{s}^{-1}$ and the red line corresponds to $k_{\text{on}}=0.023\,\text{s}^{-1}$.
It is observed that the bridge grows with a faster rate as we increase the pili-pili binding rate.
The right plot renders the bridge length versus time for three different values of the pili-pili mediated attractive forces with $k_{\text{on}}=0.023\,\text{s}^{-1}$.
The blue line corresponds to $f^{\text{p}}=18\,\text{pN}$, the green line corresponds to $f^{\text{p}}=31\,\text{pN}$ and the red line corresponds to $f^{\text{p}}=43\,\text{pN}$.
Similarly, larger values of $f^{\text{p}}$ result in a faster growth of the bridge.
These two observations are justifiable since larger binding rates and larger pili-pili attractive forces imply quicker bound pili formation and stronger attractive forces between the cells, respectively, and thus they yield quicker aggregate formation.
In both cases, the growth starts with smaller rate due to the function of the passive stress.
Further evolution of time leads the pili-pili mediated forces to become more dominant which yields higher rate growths.

\subsection{Aggregate position and periodicity}\label{sec:periodicity}

\begin{figure}[b!]
\centering
\includegraphics[width=1.0\textwidth]{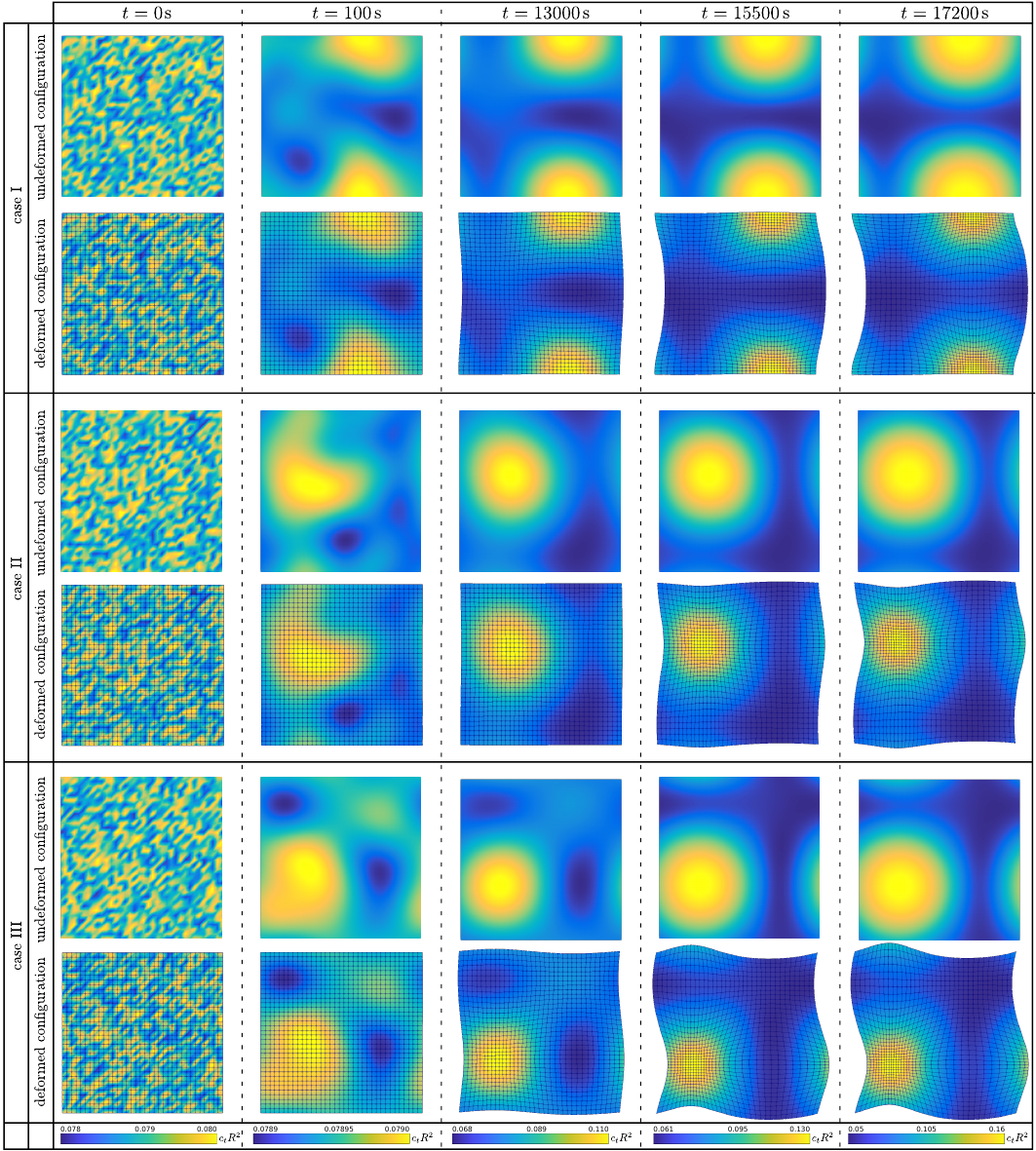}
\caption{
Formation of aggregates under different initial random distributions.
In each block, the top figures render the undeformed configuration whereas the bottom figures render the deformed bodies.
Although the initial distributions are the different, aggregates with similar sizes and same degree of phase separations are obtained. 
}
\label{fig:position}
\end{figure}

\begin{figure}[t!]
\centering
\includegraphics[width=1.0\textwidth]{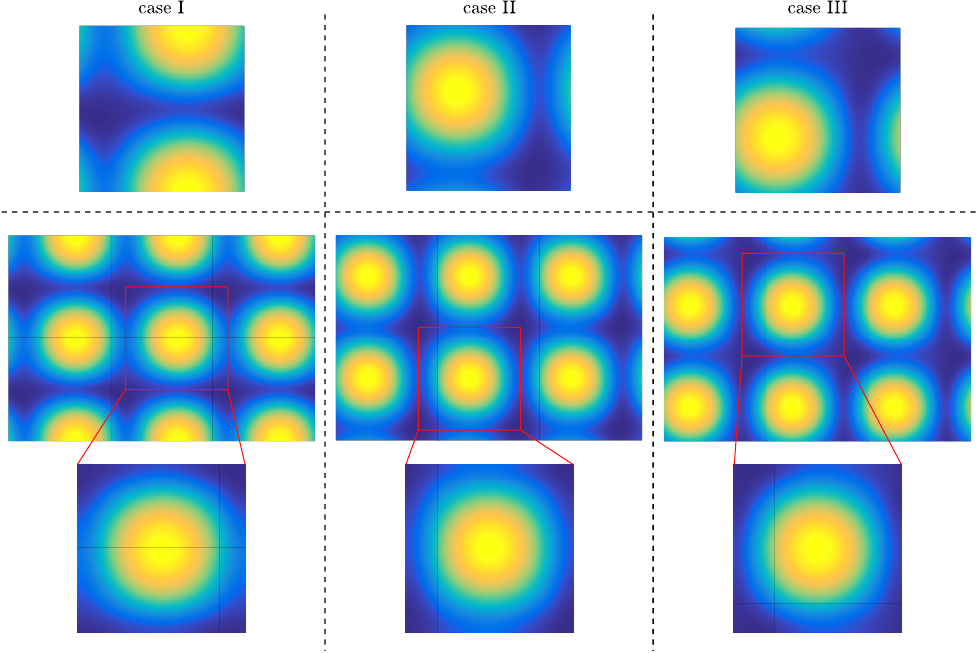}
\caption{
Illustration of the periodicity of the domain.
Although the aggregate are formed in different positions in the domain, due to periodicity, a single centered aggregate can always be extracted if enough samples are put next to each other.
}
\label{fig:periodic}
\end{figure}

\noindent
This section shows the formation of aggregates under different initial random distribution of cells.
Figure~\ref{fig:position} renders three different cases with various initial random distribution of cells.
Each block contains the undeformed and deformed configurations with five different snapshots from the domain.
It is observed that, despite different initial conditions, similar aggregates with the same size and degree of the phase separation are obtained at the end.
It is noteworthy that due to the different initial distribution, the bodies have deformed differently in order to yield the same aggregate.
Although the aggregates are formed in different locations in the three cases shown in Fig.~\ref{fig:position}, due to periodicity, a single centered aggregate can always be extracted if enough samples are put next to each other. 
Figure~\ref{fig:periodic} sheds light on this issue more vividly.
The final snapshots of the three cases in Fig.~\ref{fig:position} are shown on the top row in Fig.~\ref{fig:periodic} .
In the middle row, these samples are put together in order to form a periodic structure.
In the bottom row, it is shown that similar centered aggregates can always be extracted.


\section{Summary and outlook}\label{sec:conc}
\noindent
A continuum framework to model and simulate the behavior of biological cellular aggregates has been established.
The process of micro-colony formation has been described as an active phase separation phenomenon.
It turns out that employing the Lagrangian approach yields considerable simplification of the equations in particular for the active stress time evolution as compared to previously introduced Eulerian approach~\cite{Kuan2021}.
In addition to satisfying the conservation of the total cell number, our proposed Lagrangian formulation enabled implicit time integration which considerably increases the computational robustness.
We demonstrate that three major forces determine the dynamics of the cells in an aggregate network; the pili-pili mediated attractive forces, the steric repulsion forces and the cell-substrate friction.
In the absence of pili-pili mediated forces, the repulsive forces simply distribute the cells uniformly throughout the domain whereas in the presence of the pili-pili mediated forces, we observe a phase separation leading to the formation of micro-colonies.
A parametric study has been carried out to study the influence of various parameters on the cell density behavior.
Our proposed methodology furnishes a general framework for the continuum modeling of the non-equilibrium dynamics of dense cellular aggregates.
We believe that this contribution provides significant insights towards the dynamics of cell aggregates which in turn can be exploited to better understand the behavior of infectious diseases.
Further extension of this work include analysis of different cell species, as well as examining other biological systems with pronounced cell-matrix interactions.

\section*{Acknowledgment}
\noindent 
Soheil Firooz and Paul Steinmann gratefully acknowledge the support provided by EAM cluster.
Also Soheil Firooz would like to thank Hui--Shun Kuan for fruitful discussions regarding the parametric study.
Vasily Zaburdaev would like to acknowledge the support by Volkswagen foundation ``Life?'' initiative.

\appendix

\section{Detailed derivations of the weak forms}\label{sec:app_1}
\noindent
This section provides further details regarding the intermediate steps in the derivation of the linear momentum balance equation~\eqref{eq:weak_LMB_Lagrangian_2} as follows
\begin{equation}
\begin{aligned}
\int_{\c{B}_0} \nabla_{\b{X}}\cdot &\b{P}^{\text{a}} \cdot \delta \b{y} \, \d V + \int_{\c{B}_0} \nabla_{\b{X}} \cdot \b{P}^{\text{p}} \cdot \delta \b{y} \, \d V - \int_{\c{B}_0} \xi J c_{t} \dot{\b{y}} \cdot \delta \b{y} \, \d V\\
=
&\int_{\c{B}_0}  \nabla_{\b{X}} \cdot \lp \b{P}^{a^{T}} \cdot \delta \b{y} \rp \, \d V - \int_{\c{B}_0}  \b{P}^{\text{a}} : \nabla_{\b{X}} \delta \b{y}  \, \d V + \int_{\c{B}_0}  \nabla_{\b{X}} \cdot \lp \b{P}^{p^{T}} \cdot \delta \b{y} \rp \, \d V - \int_{\c{B}_0}  \b{P}^{\text{p}} : \nabla_{\b{X}} \delta \b{y}  \, \d V  - \int_{\c{B}_0} \xi J c_{t} \dot{\b{y}} \cdot \delta \b{y} \, \d V\\
=
&\int_{\c{B}_0}  \lb \b{P}^{a^{T}} \cdot \delta \b{y} \rb \cdot \b{N} \, \d A - \int_{\c{B}_0}  \b{P}^{\text{a}} : \nabla_{\b{X}} \delta \b{y}  \, \d V + \int_{\c{B}_0}  \lb \b{P}^{p^{T}} \cdot \delta \b{y} \rb \cdot \b{N} \, \d A - \int_{\c{B}_0}  \b{P}^{\text{p}} : \nabla_{\b{X}} \delta \b{y}  \, \d V  - \int_{\c{B}_0} \xi J c_{t} \dot{\b{y}} \cdot \delta \b{y} \, \d V\\
=
&\int_{\c{B}_0}   \delta \b{y} \cdot \b{P}^{\text{a}} \cdot \b{N} \, \d A - \int_{\c{B}_0}  \b{P}^{\text{a}} : \nabla_{\b{X}} \delta \b{y}  \, \d V + \int_{\c{B}_0}   \delta \b{y} \cdot \b{P}^{\text{p}} \cdot \b{N} \, \d A - \int_{\c{B}_0}  \b{P}^{\text{p}} : \nabla_{\b{X}} \delta \b{y}  \, \d V  - \int_{\c{B}_0} \xi J c_{t} \dot{\b{y}} \cdot \delta \b{y} \, \d V\\
=
&\int_{\c{B}_0}   \delta \b{y} \cdot \b{T}^{\text{a}}\, \d A - \int_{\c{B}_0}  \b{P}^{\text{a}} : \nabla_{\b{X}} \delta \b{y}  \, \d V + \int_{\c{B}_0}   \delta \b{y} \cdot \b{T}^{\text{p}} \, \d A - \int_{\c{B}_0}  \b{P}^{\text{p}} : \nabla_{\b{X}} \delta \b{y}  \, \d V  - \int_{\c{B}_0} \xi J c_{t} \dot{\b{y}} \cdot \delta \b{y} \, \d V = 0\,.
\end{aligned}
\label{eq:weak_LMB_Lagrangian_1_app}
\end{equation}

\section{Detailed derivations of the active stress}\label{sec:app_2}
\noindent
In this section, the derivation of the material time derivative of the active Piola--Kirchhoff stress is provided in detail.
To begin, the Lie time derivative of the Kirchhoff stress is obtained as follows
\begin{equation}
\begin{aligned}
\c{L}_t\b{\tau}^\text{a} 
&=  \b{F} \cdot \lb \dot{\o{\b{F}^{-1} \cdot\b{\tau}^\text{a}\cdot \b{F}^{-T}}} \rb  \cdot \b{F}^T
 = \b{F} \cdot \left[  \dot{\b{F}^{-1}} \cdot \b{\tau}^\text{a} \cdot \b{F}^{-T} + \b{F}^{-1}  \cdot \dot{\b{\tau}^\text{a}} \cdot \b{F}^{-T} + \b{F}^{-1} \cdot \b{\tau}^\text{a} \cdot \dot{\b{F}^{-T}}  \right] \cdot \b{F}^T\\[8pt]
& = \b{F}\cdot \left[-  \b{F}^{-1}\cdot\b{l}\cdot \b{\tau}^\text{a}\cdot \b{F}^{-T} + \b{F}^{-1}\cdot \dot{\b{\tau}^\text{a}} \cdot \b{F}^{-T} - \b{F}^{-1} \cdot \b{\tau}^\text{a} \cdot \b{l}^{T} \cdot \b{F}^{-T}   \right] \cdot \b{F}^T \\[8pt]
& = \dot{\b{\tau}^\text{a}} - \b{l} \cdot \b{\tau}^\text{a}  - \b{\tau}^\text{a} \cdot \b{l}^{T} 
  = \dot{\b{\tau}^\text{a}} - \b{l}\cdot {\b{\tau}^\text{a}}^T  - \b{\tau}^\text{a}\cdot \b{l}^{T} 
  = \dot{\b{\tau}^\text{a}} - 2\left[ \b{l} \cdot \b{\tau}^\text{a} \right]^{\text{sym}}\,.
\end{aligned}
\label{eq:Lie_time_derivative_app}
\end{equation}
Afterwards, using the relation $\b{S}=\b{F}^{-1} \cdot \b{\tau}^\text{a} \cdot \b{F}^{-T}$ on could write
\begin{equation}
\begin{aligned}
\dot{\b{S}^\text{a}}
= - \f{1}{\ell_0 p_0f^{\text{p}}} \left[ \lb\b{F}\cdot\b{S}^\text{a}\cdot\b{F}^{T}\rb:\b{l}^{T}\right]\b{S}^\text{a}  + \b{S}^\text{f} - k_\text{off}\,\b{S}^\text{a}
= - \f{1}{\ell_0 p_0f^{\text{p}}} \left[ \b{S}^\text{a}:\lb\b{F}^{T}\cdot\b{l}^{T}\cdot\b{F}\rb\right]\b{S}^\text{a}  + \b{S}^\text{f} - k_\text{off}\,\b{S}^\text{a}\,.
\end{aligned}
\label{eq:active_stress_Lagrangian_app_1}
\end{equation}
Since the Piola--Kirchhoff stress is a symmetric tensor, for an arbitrary second order tensor $\b{A}$ we can write
\begin{equation}
\b{S}:\b{A} = \b{S}:\left[\f{1}{2}\left[ \b{A} + \b{A}^T \right]\right]\,.
\end{equation}
Thus, Eq.~\eqref{eq:active_stress_Lagrangian_app_1} can be rewritten as
\begin{equation}
\begin{aligned}
\dot{\b{S}^\text{a}}
&= - \f{1}{\ell_0 p_0f^{\text{p}}} \left[ \b{S}^\text{a}:\lb\b{F}^{T}\cdot\b{l}^{T}\cdot\b{F}\rb\right]\b{S}^\text{a}  + \b{S}^\text{f} - k_\text{off}\,\b{S}^\text{a}\\
&= - \f{1}{\ell_0 p_0f^{\text{p}}} \left[ \b{S}^\text{a}:\f{1}{2}\left[\lb\b{F}^{T}\cdot\b{l}^{T}\cdot\b{F}\rb+ \left[\b{F}^{T}\cdot\b{l}^{T}\cdot\b{F}\right]^{T}\right]\right]\b{S}^\text{a}  + \b{S}^\text{f} - k_\text{off}\,\b{S}^\text{a}\\
&= - \f{1}{\ell_0 p_0f^{\text{p}}} \left[ \b{S}^\text{a}:\lb\b{F}^{T}\cdot\b{l}^{\text{sym}}\cdot\b{F}\rb\right]\b{S}^\text{a}  + \b{S}^\text{f} - k_\text{off}\,\b{S}^\text{a}\\
\end{aligned}
\label{eq:active_stress_Lagrangian_app_2}
\end{equation}
Finally, utilizing the relation
\begin{equation}
\b{F}^{T}\cdot[\b{l}]^\text{sym}\cdot\b{F} = 
\b{F}^T \cdot\left[\f{1}{2}\left[\b{l}^T + \b{l} \right] \right]\cdot\b{F} =
\f{1}{2}\left[\b{F}^T\cdot\b{l}^T\cdot\b{F} + \b{F}^T\cdot\b{l}\cdot\b{F} \right] =
\f{1}{2}\left[\dot{\b{F}^T}\cdot\b{F} + \b{F}^T\cdot\dot{\b{F}} \right] =
\dot{\b{E}}\,,
\end{equation}
we can derive the fully Lagrangian form of the active Piola--Kirchhoff stress as
\begin{equation}
\dot{\b{S}^\text{a}}
= - \f{1}{\ell_0 p_0f^{\text{p}}} \left[ \b{S}^\text{a}:\dot{\b{E}}\right]\b{S}^\text{a}  + \b{S}^\text{f} - k_\text{off}\,\b{S}^\text{a}\,.
\label{eq:active_stress_Lagrangian_app_3}
\end{equation}

\section{Time integration of the active stress}\label{sec:app_3}
\noindent
In this section we detail on the time integration technique in order to calculate the active second Piola--Kirchhoff stress.
The non-linear relation for the active second Piola--Kirchhoff stress reads
\begin{equation}
\dot{\b{S}^\text{a}} 
= - \f{1}{\ell_0 p_0f^{\text{p}}} \left[ \b{S}^\text{a}:\dot{\b{E}}\right]\b{S}^\text{a}  + \b{S}^\text{f} - k_\text{off}\,\b{S}^\text{a}\,,
\end{equation}
which could be written in the form
\begin{equation}
\f{\b{S}^{\text{a}}_{n+1}-\b{S}^{\text{a}}_{n}}{\Delta t} = - \f{1}{\ell_{0} p_{0}f^{\text{p}}\Delta t} \b{S}^{\text{a}}_{n+1} \lb \b{S}^{\text{a}}_{n+1}:\lb \b{E}_{n+1} - \b{E}_{n} \rb \rb + \b{S}_{n+1}^{\text{f}} - k_{\text{off}}\b{S}^{\text{a}}_{n+1}\,.
\end{equation}
To linearize this equation, we put all the terms on one side and treat them as a residuum $\b{R}$ that must vanish
\begin{equation}
\b{R} = \f{\b{S}^{\text{a}}_{n+1}-\b{S}^{\text{a}}_{n}}{\Delta t} + \f{1}{\ell_{0} p_{0}f^{\text{p}}\Delta t} \b{S}^{\text{a}}_{n+1} \lb \b{S}^{\text{a}}_{n+1}:\lb \b{E}_{n+1} - \b{E}_{n} \rb \rb - \b{S}_{n+1}^{\text{f}} + k_{\text{off}}\b{S}^{\text{a}}_{n+1}\,.
\end{equation}
The linearization of $\b{R}$ reads
\begin{equation}
\b{R}_{n+1} \approx \text{Lin}\b{R}_{n+1} = \b{R}_{n} + \f{\p\b{R}}{\p\b{S}^{\text{a}}}\big\lvert_{n} \cdot \Delta \b{S}^{\text{a}}_{n}\stackrel{!}{=}\bnull\,.
\label{eq:linearize_R}
\end{equation}
The tangent reads
\begin{equation}
\begin{aligned}
\v{K} = \f{\p\b{R}}{\p\b{S}^{\text{a}}} &= \lb \f{1}{\Delta t} + k_{\text{off}}\rb \v{I}^{\text{sym}} + \f{1}{\ell_{0} p_{0}f^{\text{p}}\Delta t}  \bigg[ \lb \b{S}^{\text{a}}_{n+1}:\lb \b{E}_{n+1} - \b{E}_{n} \rb \rb  \v{I}^{\text{sym}} + \b{S}^{\text{a}}_{n+1}\dyad\lb \b{E}_{n+1} - \b{E}_{n} \rb \bigg]\\
&=\lb \f{1}{\Delta t} + k_{\text{off}}  + \f{1}{\ell_{0} p_{0}f^{\text{p}}\Delta t} \lb \b{S}^{\text{a}}_{n+1}:\lb \b{E}_{n+1} - \b{E}_{n} \rb \rb \rb \v{I}^{\text{sym}} + \f{1}{\ell_{0}p_{0}f^{\text{p}}\Delta t}  \bigg[   \b{S}^{\text{a}}_{n+1}\dyad\lb \b{E}_{n+1} - \b{E}_{n} \rb \bigg]\,,
\end{aligned}
\end{equation}
with $\v{I}^{\text{sym}}$ being the symmetric fourth-order identity which reads $\v{I}^{\text{sym}}=1/2\lb \delta_{ik}\delta_{jl} + \delta_{il}\delta_{jk} \rb$.
Finally, an iterative Newton--Raphson scheme is employed to solve for $\b{R}_{n+1}$ and thus, for $\b{S}^{\text{a}}_{n+1}$ as
\begin{equation}
\lb\b{S}^{\text{a}}\rb_{n+1} = \lb\b {S}^{\text{a}} \rb_{n} + \lb d\b{S}^{\text{a}} \rb \quad \text{with} \quad  \lb \Delta\b{S}^{\text{a}} \rb_{ij} = - \lb \v{K}^{-1} \rb_{ijkl} \lb \b{R} \rb_{kl}\,.
\end{equation}
To calculate $\v{K}^{-1}$ we use the Sherman--Morrison formula which states that for an arbitrary fourth-order tensor $\v{A}$ that can be written in the form 
\begin{equation}
\v{A} = \beta \v{B} + \alpha \b{C}\dyad\b{D} \,,
\label{eq:inverse}
\end{equation}
with $\alpha$ and $\beta$ being scalars, $\b{C}$ and $\b{D}$ being second-order tensors and $\v{B}$ being a fourth-order tensor, the inverse of $\v{A}$ reads
\begin{equation}
\v{A}^{-1} = \f{1}{\beta} \v{B}^{-1} -  \f{\alpha}{\beta^2+\alpha\beta\b{D}:\v{B}^{-1}:\b{C}} \lb \v{B}^{-1}:\b{C}\dyad\b{D}:\v{B}^{-1} \rb\,.
\end{equation}
To proceed, we write our tangent in the form of Eq.~\eqref{eq:inverse} as 
\begin{equation}
\v{K} = \beta \v{I}^{\text{sym}} + \alpha \bigg[   \b{S}_{n+1}\dyad\lb \b{E}_{n+1} - \b{E}_{n} \rb \bigg]\,,
\end{equation}
with 
\begin{equation}
\alpha=\f{1}{\ell_{0}\rho^{\text{p}}_{0}f^{\text{p}}\Delta t} \qquad \text{and} \qquad 
\beta=\lb \f{1}{\Delta t} + k_{\text{off}}  + \f{1}{\ell_{0}\rho^{\text{p}}_{0}f^{\text{p}}\Delta t} \lb \b{S}_{n+1}:\lb \b{E}_{n+1} - \b{E}_{n} \rb \rb \rb\,.
\end{equation}
Therefore, using the identity $\v{I}^{\text{sym}^{-1}}=\v{I}^{\text{sym}}$, $\v{K}^{-1}$ reads
\begin{equation}
\v{K}^{-1} = \f{1}{\beta} \v{I}^{\text{sym}} -  \f{\alpha}{\beta^2+\alpha\beta\lb \b{E}_{n+1} - \b{E}_{n} \rb:\v{I}^{\text{sym}}:\b{S}_{n+1}} \lb \v{I}^{\text{sym}}:\b{S}_{n+1}\dyad\lb \b{E}_{n+1} - \b{E}_{n} \rb:\v{I}^{\text{sym}} \rb\,,
\end{equation}
or in index notation
\begin{equation}
\lb\v{K}^{-1}\rb_{ijkl} \!\!\!= \f{1}{\beta}\lb \v{I}^{\text{sym}}\rb_{ijkl} -  \f{\alpha}{\beta^2+\alpha\beta\lb \b{E}_{n+1} - \b{E}_{n} \rb_{mn}:\lb\v{I}^{\text{sym}}\rb_{mnrs}:\lb\b{S}_{n+1}\rb_{rs}  } \lb \lb\v{I}^{\text{sym}}\rb_{ijmn}: \lb\b{S}_{n+1}\dyad\lb \b{E}_{n+1} - \b{E}_{n} \rb\rb_{mnrs}:\lb\v{I}^{\text{sym}}\rb_{rskl} \rb\,.
\end{equation}

\section{Some useful relations}\label{sec:app_4}
\noindent
This section provides some useful derivatives in index notation relation which will be helpful regarding the derivation of the tangents for the finite element implementation.
\begin{equation*}
\left[\f{\p J}{\p \b{F}}\right]_{ij} = J \left[\b{F}^{-T}\right]_{ij} \,,
\qquad
\left[\f{\p \b{F}}{\p \b{F}}\right]_{ijkl} = \delta_{ik} \delta_{jl}\,,
\qquad
\left[\f{\p \b{F}^{T}}{\p \b{F}}\right]_{ijkl} = \delta_{il} \delta_{jk}\,,
\qquad
\left[\f{\p \b{F}^{-1}}{\p \b{F}}\right]_{ijkl} = - \left[\b{F}^{-1} \right]_{ik} \left[\b{F}^{-T} \right]_{jl}\,,
\end{equation*}
\vspace{-0.5cm}
\begin{equation*}
\left[\f{\p \b{F}^{-T}}{\p \b{F}}\right]_{ijkl} = - \left[\b{F}^{-1} \right]_{il} \left[\b{F}^{-T} \right]_{jk}\,,
\qquad
\left[\f{\p \b{B}}{\p \b{F}}\right]_{ijkl} = \left[\f{\p \b{F}^{-1}}{\p \b{F}}\right]_{imkl} \left[\b{F}^{-T} \right]_{mj} + \left[\b{F}^{-1} \right]_{im} \left[\f{\p \b{F}^{-T}}{\p \b{F}}\right]_{mjkl} \,,
\end{equation*}
\vspace{-0.5cm}
\begin{equation*}
\left[\f{\p \b{E}}{\p \b{F}}\right]_{ijkl} = \f{1}{2} \left[ \left[\f{\p \b{F}^{T}}{\p \b{F}}\right]_{imkl} \left[\b{F} \right]_{mj} + \left[\b{F}^{T} \right]_{im} \left[\f{\p \b{F}}{\p \b{F}}\right]_{mjkl} \right]\,,\\
\end{equation*}
\vspace{-0.5cm}
\begin{equation*}
\left[\f{\p \b{P}}{\p \b{F}}\right]_{ijkl} = \left[\f{\p \b{F}}{\p \b{F}}\right]_{imkl} \left[\b{S} \right]_{mj} + \left[\b{F} \right]_{im} \left[\f{\p \b{S}}{\p \b{F}}\right]_{mjkl} \,,
\qquad
\left[ \f{\p [\b{g} \dyad \b{g}]}{\p \b{g}} \right]_{ijk} \!\!\!\!= \left[ \left[\b{I}\right]_{ik} \left[\b{g}\right]_{j}  + \left[\b{g}\right]_{i} \left[\b{I}\right]_{jk} \right]\,,
\end{equation*}
\vspace{-0.5cm}
\begin{equation*}
\left[ \f{\p [\nabla_{\b{X}}\,\b{g} : \b{F}^{-T}]}{\p \nabla_{\b{X}}\,\b{g}} \right]_{ij} \!\!\!\!= \left[\b{F}^{-T}\right]_{ij} \left[ \f{\p [\nabla_{\b{X}}\,\b{g} : \b{F}^{-T}]}{\p\,\b{F}} \right]_{ij} \!\!\!\!=  \left[\f{\p\b{F}^{-T}}{\p\b{F}}\right]_{klij} \left[\nabla_{\b{X}}\,\b{g} \right]_{kl}\,.
\end{equation*}

\section{Calculation of the derivatives at the element level}\label{sec:app_5}
\noindent
In this section a brief example to calculate the derivatives of a scalar field $\alpha$, a vector field $\b{a}$ and a second-order tensor field $\b{A}$ with respect to the nodal values at the element level is elaborated.
\begin{equation*}
\lb \f{\p \alpha}{\p \alpha^{J}} \rb = \f{\p\lp \alpha^{S}N^{S}\rp}{\p \alpha^{J}} = \f{\p \alpha^{S}}{\p \alpha^{J}} N^{S} =  \delta^{SJ} N^{S} = N^{J}\,,
\end{equation*}
\vspace{-0.5cm}
\begin{equation*}
\lb \f{\p\nabla_{\b{x}} \alpha}{\p \alpha^{J}} \rb_{i} = \f{\p\lb\nabla_{\b{x}}\lp \alpha^{S}N^{S}\rp\rb_{i}}{\p \alpha^{J}} = \f{\p\lb \alpha^{S} \nabla_{\b{x}} N^{S}_{i} \rb}{\p \alpha^{J}} = \f{\p \alpha^{S}}{\p \alpha^{J}} \, \nabla_{\b{x}} N^{S}_{i} = \delta^{SJ} \, \nabla_{\b{x}} N^{S}_{i} = \lb \nabla_{\b{x}} N^{J}\rb_{i}\,,
\end{equation*}
\vspace{-0.5cm}
\begin{equation*}
\lb \f{\p\b{a}}{\p\b{a}^{J}} \rb_{ij} = \f{\p\lp\b{a}^{S}_{i}N^{S}\rp}{\p\b{a}^{J}_{j}} = \f{\p\b{a}^{S}_{i}}{\p\b{a}^{J}_{j}} N^{S} = \delta_{ij} \delta^{SJ} N^{S} = \delta_{ij} N^{J}\,,
\end{equation*}
\vspace{-0.5cm}
\begin{equation*}
\lb \f{\p\nabla_{\b{x}}\b{a}}{\p\b{a}^{J}} \rb_{ijk} = \f{\p\lb\nabla_{\b{x}}\lp\b{a}^{S}_{i}N^{S}\rp\rb_{j}}{\p\b{a}^{J}_{k}} = \f{\p\lb\b{a}^{S}_{i} \dyad \nabla_{\b{x}} N^{S}_{j} \rb}{\p\b{a}^{J}_{k}} = \f{\p\b{a}^{S}_{i}}{\p\b{a}^{J}_{k}} \, \nabla_{\b{x}} N^{S}_{j} = \delta_{ik} \delta^{SJ} \, \nabla_{\b{x}} N^{S}_{j} = \delta_{ik} \, \lb\nabla_{\b{x}} N^{J}\rb_{j}\,,
\end{equation*}
\vspace{-0.5cm}
\begin{equation*}
\lb \f{\p\b{A}}{\p\b{A}^{J}} \rb_{ijkl}\!\!\!\! = \f{\p\lp\b{A}^{S}_{ij}N^{S}\rp}{\p\b{A}^{J}_{kl}} = \f{\p\b{A}^{S}_{ij}}{\p\b{A}^{J}_{kl}} N^{S} = \delta_{ik} \delta_{jl} \delta^{SJ} N^{S} = \delta_{ik} \delta_{jl} N^{J}\,,
\end{equation*}
\vspace{-0.5cm}
\begin{equation*}
\lb \f{\p\nabla_{\b{x}}\b{A}}{\p\b{A}^{J}} \rb_{ijklm} = \f{\p\lb\nabla_{\b{x}}\lp\b{A}^{S}_{ij}N^{S}\rp\rb_{k}}{\p\b{A}^{J}_{lm}} = \f{\p\lb\b{A}^{S}_{ij} \dyad \nabla_{\b{x}} N^{S}_{k} \rb}{\p\b{A}^{J}_{lm}} = \f{\p\b{A}^{S}_{ij}}{\p\b{A}^{J}_{lm}} \, \nabla_{\b{x}} N^{S}_{k} = \delta_{il} \delta_{jm} \delta^{SJ} \, \nabla_{\b{x}} N^{S}_{k} = \delta_{il} \delta_{jm} \, \lb\nabla_{\b{x}} N^{J}\rb_{k}\,.
\end{equation*}

\section*{References}
\bibliography{library}

\end{document}